\def\WHO{nobody}
\def\version{19.7.2024}
\def\users{us}  %
\def\users{final-layout}   % when activated, ``our'' debugging is suppressed
\documentclass[12pt]{article}
\usepackage{upgreek}
\usepackage{psfrag}    % for figures
\usepackage{graphicx}  % standard LaTeX graphics tool
\usepackage{xcolor}
\usepackage{bm,amsmath,amsthm,hyperref,amsfonts,amssymb,color}
\usepackage{mathrsfs} % THIS PACKAGE ALLOWS FOR ``scr'' FONTS
\usepackage{cite}

\usepackage{ifthen}
\ifthenelse{\equal{\users}{final-layout}}{}{
\usepackage{fancyhdr}
\pagestyle{fancy}
\headheight=28pt\headwidth=17cm
\definecolor{gray}{gray}{0.5}
\rhead{\color{gray}Planetary accretion \\
T.Roub\'\i\v cek \& U.Stefanelli}
\chead{}
\lhead{Version \version, file: \jobname.tex,
\\
compiled:
\number\day.\number\month.\number\year\ at
\the\hour:\ifnum\minute<10 0\fi\the\minute\ h\ \ \ \ \ 
\underline{\Large\color{red}\WHO's working on!}
}
}

\definecolor{labelkey}{rgb}{1.,.2,0.}

\newcount\hour \newcount\minute
\hour=\time
\divide \hour by 60
\minute=\time
\loop \ifnum \minute > 59 \advance \minute by -60 \repeat

\usepackage[normalem]{ulem}

\usepackage[normalem]{ulem}
\usepackage{ifthen}
\usepackage{color}

\ifthenelse{\equal{\users}{final-layout}}{

	\newcommand{\COMMENT}[1]{}
	\newcommand{\COMMENTGT}[1]{}
	\newcommand{\TODO}[1]{}
	\newcommand{\INTERNAL}[1]{}
	\newcommand{\QUESTION}[1]{}
	\newcommand{\DELETE}[1]{}

	\newcommand{\REM}[1]{\marginpar{\bfseries\tiny{\color{blue}}}}
    \newcommand{\MARGINOTE}[1]{}
}
{
	
	\newcommand{\COMMENT}[1]{{\color{red}\uuline{#1}\color{black}}}
	\newcommand{\COMMENTGT}[1]{{\hfill\large\color{red}***{#1}***\color{black}\hfill}\\}
	\newcommand{\TODO}[1]{{\color{red}\uuline{#1}\color{black}}}
	\newcommand{\INTERNAL}[1]{\footnote{#1}}
	\newcommand{\QUESTION}[1]{{\color{brown}\uuline{#1}\color{black}}}
	\newcommand{\DELETE}[1]{{\color{red}\sout{#1}\color{black}}}

	\newcommand{\REM}[1]{\marginpar{\bfseries\tiny{\color{blue}#1}}}
\newcommand{\MARGINOTE}[1]{\marginpar{\color{red}\tiny\texttt{#1}}}
}

\newcommand\DT[1]{\mathchoice
                 {{\buildrel{\hspace*{.1em}\text{\LARGE.}}\over{#1}}}
                 {{\buildrel{\hspace*{.1em}\text{\LARGE.}}\over{#1}}}
                 {{\buildrel{\hspace*{.1em}\text{\Large.}}\over{#1}}}
                 {{\buildrel{\hspace*{.1em}\text{\large.}}\over{#1}}}}
\newcommand\pdt[1]{\frac{\partial{#1}}{\partial t}} %Partial Derivative w.r.t. t
\newcommand{\lineunder}[2]{\LU{\begin{array}[t]{c}\underbrace{#1}\vspace*{.5em}\end{array}}{\mbox{\footnotesize\rm #2}}}
\newcommand{\LU}[2]{\begin{array}[t]{c}#1\vspace*{-1em}\\_{#2}\end{array}}
\newcommand{\linesunder}[3]{\LSU{\begin{array}[t]{c}\underbrace{#1}\vspace*{.5em}\end{array}}{\mbox{\footnotesize\rm #2}}{\mbox{\footnotesize\rm #3}}}
\newcommand{\LSU}[3]{\begin{array}[t]{c}#1\vspace*{-1em}\\_{#2}\vspace*{-.5em}\\_{#3}\end{array}}
\newcommand{\morelinesunder}[4]{\LSUU{\begin{array}[t]{c}\underbrace{#1}\vspace*{.5em}\end{array}}{\mbox{\footnotesize\rm #2}}{\mbox{\footnotesize\rm #3}}{\mbox{\footnotesize\rm #4}}}
\newcommand{\LSUU}[4]{\begin{array}[t]{c}#1\vspace*{-1em}\\_{#2}\vspace*{-.5em}\\_{#3}\vspace*{-.5em}\\_{#4}\end{array}}

\def\vv{\bm v}
\def\oomega{\bm \omega}
\def\XX{\bm X}
\def\xx{\bm x}
\def\vv{\bm v}

\def\gg{\bm g}
\def\nn{\bm n}
\def\ff{\bm f}
\def\jj{\bm j}
\def\yy{\bm y}
\def\qq{\bm q}
\def\ee{\bm e}
\def\FF{\bm F}
\def\DD{\bm D}
\def\TT{\bm T}
\def\GRAVCONST{G}
\def\GRAVPOT{V}
\def\bbI{ \bm I}
\newcommand{\eq}[1]{(\ref{#1})}
\newcommand{\Cdot}{\hspace{-.1em}\cdot\hspace{-.1em}}

\def\d{{\rm d}}
\def\W{w}
\def\COUPLING{\gamma}
\def\DIS{{\large\mbox{$\delta$}}}
\def\nothing{}

\usepackage{mathrsfs}   % loading \mathscr
\usepackage{eucal}      % changing \mathcal into EulerCaligraphic  
\newcommand{\R}{\mathbb R}

\newcommand{\haz}{\widehat}

\newtheorem{theorem}{Theorem}[section]

\newtheorem{example}[theorem]{Example}

\newtheorem{remark}[theorem]{Remark}

\newcommand{\ttt}{\color{black}}
\newcommand{\sss}{\color{black}}
\newcommand{\eee}{\color{black}}

\numberwithin{equation}{section}

\begin{document}

\allowdisplaybreaks

\bigskip\bigskip\bigskip

\noindent{\LARGE\bf
Modelling of planetary accretion and\\[.2em]core-mantle structure formation
}

\bigskip\bigskip\bigskip\bigskip

\noindent{\bf Tom\'a\v s Roub\'\i\v cek}\footnote{Mathematical Institute, Charles University,
Sokolovsk\'a 83, CZ-186~75~Praha~8,  Czech Republic,\\\hspace*{1.6em} email: ${\texttt{tomas.roubicek@mff.cuni.cz}}$}
\footnote{Institute of Thermomechanics, Czech Academy of Sciences,
Dolej\v skova 5, CZ-18200~Praha~8, Czech Rep.}
\& {\bf Ulisse Stefanelli}\footnote{Faculty of Mathematics, University of
  Vienna, Oskar-Morgenstern-Platz 1, 1090 Vienna, Austria and Vienna Research Platform on Accelerating Photoreaction Discovery, University of Vienna, W\"ahringerstrasse 17, A-1090 Vienna, Austria\\\hspace*{1.6em} email: ${\texttt{ulisse.stefanelli@univie.ac.at}}$}
\footnote{Istituto di Matematica Applicata e Tecnologie Informatiche
  {\it E. Magenes} - CNR, v. Ferrata 1, I-27100 Pavia,
  Italy.}\footnote{Corresponding Author.}

%\end{center}

\bigskip\bigskip\bigskip

{\small

  \noindent{\bf Abstract}.
   We advance a  thermodynamically consistent model of 
  self-gravitational accretion and differentiation
  in planets. The system is modeled in actual variables as a compressible
  thermoviscoelastic fluid in  a fixed, sufficiently large
  domain.  The supply of material to the accreting and
  differentiating system is described as a bulk source of mass,
  volume, impulse, and energy localized in some border
  region of the domain. Mass, momentum, and energy conservation, along
  with constitutive relations, result in an extended compressible
  Navier-Stokes-Fourier-Poisson system. \ttt The centrifugal and Coriolis
  forces are also considered. \eee
  After studying some single-component setting, we consider a
  two-component situation, where metals and silicates mix and
  differentiate under gravity, eventually forming a
  core-mantle structure. The energetics of the models are
  elucidated. Moreover, we prove that the models are stable, in that
  self-gravitational collapse is excluded. Eventually, we comment on
  the prospects of devising a rigorous mathematical approximation and
  existence theory.

\medskip

\medskip

\noindent {\it AMS Subject Classification:} 
35Q49, % Transport equations
35Q74, % PDEs in connection with mechanics of deformable solids
%35Q79, % PDEs in connection with classical thermodynamics and heat transfer
35Q86, % PDEs in connection with geophysics
%65M60, % PDE, IVP, ... Galerkin methods, 
%74A15, % Thermodynamics
%74A30, % Nonsimple materials
% 74A45, % Theories of fracture and damage
%74C15, %Large-strain, rate-independent theories(including nonlinear plasticity)
74Dxx, %Materials of strain-rate type.
%74F10, % Fluid-solid interactions (...aero- hydro-elasticity, porosity, etc.)
%74J30, % Nonlinear waves
%74L05, % Geophysical solid mechanics
%74R20, % Anelastic fracture and damage
%76S05. % Flows in porous media; filtration; seepage
76T06, % Liquid-liquid two component flows
80A17, %Thermodynamics of continua 
%80A20, % Heat and mass transfer, heat flow - NOT IN 2020 AMS Class.
%86A17. % Global dynamics, earthquake problems
86A99. % Geophysics None of the above, but in this section

}

\bigskip\bigskip

\baselineskip=16pt

\section{Introduction}
%        ~~~~~~~~~~~~
Modelling  the formation of planets from the protoplanetary disk 
is  a major  challenge
in theoretical and computational geophysics.  In early stages (up to
4.6 billions of years ago for the Solar System),
planetary accretion mainly results from the gas-solid (or rather
gas-fluid over such long time scales) compactification transition from the
extraplanetary dust to the planetary body. In a later stage (from 4 billions
years ago for the Solar System), accretion continues mainly due to the
impact of material at the planetary surface. The different
materials composing the planet differentiate by their weight and a
core-mantle structure is formed. All these phenomena are driven by
self-gravitation, and partially influenced by the presence of the corresponding star. 

 To model such a complex process, one would be asked to consider a
very large region, possibly including the whole stellar system. In an
effort of presenting a reduced model, with the prospect of being computationally
amenable, we reduce ourselves to consider a relatively small and
artificially chosen region $\varOmega$, orbiting around the star,
where planetary accretion and differentiation can be assumed
to take place. The material influx to this
region from the external space  
 is modeled as  
a continuous influx of material.  This is  %. This is
to be interpreted as some
time/space averaging of approaching asteroids and comets,  by 
ignoring  the  probabilistic character of this influx.

By self-gravitation,  all material in $\varOmega$ 
falls towards the barycenter of the system.
Depending on density, the system experiences a sharp
transition between a gaseous phase, corresponding to extraplanetary
dust and incoming material, and a fluid phase, describing the planet
body instead. The sharp transition region is to be identified with the planet
surface.  This is combined with the subsequent gravitational differentiation
of the planet,  where the different constituents of the planet matter
are subject to different gravity forces due to their relative weight
and hence flow within the planetary body. The  basic scenario in the
geophysics of terrestrial-type planets  is that of considering the
formation of a metallic core and a silicate mantle by
self-gravitating differentiation,  see,
 e.g., \cite{GerYue07RCMM,GGKZ09RCTC}. A major role is played here by
heat, which is originating from the impact of the material on the
surface, as well as from the friction of the flowing planet constituents
during differentiation. Kinetic energy is also part of the picture and
expected to generate planetary rotation. Finally, the
orbiting of the system around the orbital center generates tidal
effects, which are here taken into account, up to some simplification,
see Subsection \ref{sec:gravity}.

In addition, the
differentiating body is subject to chemical reactions, radiation, and
electromagnetic effects. \sss These are all argued to possibly have
an impact on the planet-formation process. Referring to the Earth,
the magnetic field generated by convection of iron and nickel in the
outer core and its interaction with the solar wind are indeed
influencing thermal evolution 
\cite{McDYos21TPCC,SLKA15MDCP}.  
Aiming at maximal simplicity, all these
effects are currently neglected in this paper. %We however plan to
%include some of these effects in some future version of the model.
% We , which are all neglected in this paper.
% \sss Even during the accretion in early stages, the (here neglected) magnetic
% field may play some role due to the magnetic field comming from the central
% star and later the own magnetic field of the planet caused and backward
% influcening also the thermal evolution, cf.\ e.g.\ \cite{McDYos21TPCC,SLKA15MDCP}. 
\eee

\begin{figure}[ht]
\begin{center}
%\begin{my-picture}{6}{6.5}{fig1}
\psfrag{W}{\large$\varOmega$}
\psfrag{G}{\large$\varGamma$}
\psfrag{viz Remark 2.2}{\footnotesize (viz Remark~2.3)}
\psfrag{border}{\footnotesize border}
\psfrag{zone}{\footnotesize zone}
\psfrag{V=0}{}\psfrag{v=0}{\small $\vv{\cdot}\nn=0$} 
\psfrag{v,J}{\small $\varrho,\,\vv,\,J,\,\theta$}
\psfrag{r-ext}{\small $v_{\rm ext}^{},\,\vv_{\rm ext}^{},\,h_{\rm ext}^{}$}
\psfrag{flux=0}{\small $\nn{\dot}\nabla\theta=0$}
\psfrag{planet}{\footnotesize\begin{minipage}[t]{11em}accreting and\\[-.2em]\hspace*{-.2em}differentiating\\[-.2em]\hspace*{1.5em}planet\end{minipage}}
%\hspace*{-.5em}\includegraphics[width=17em]{tr-us-accretion1.eps}
\hspace*{-.5em}\includegraphics[width=17em]{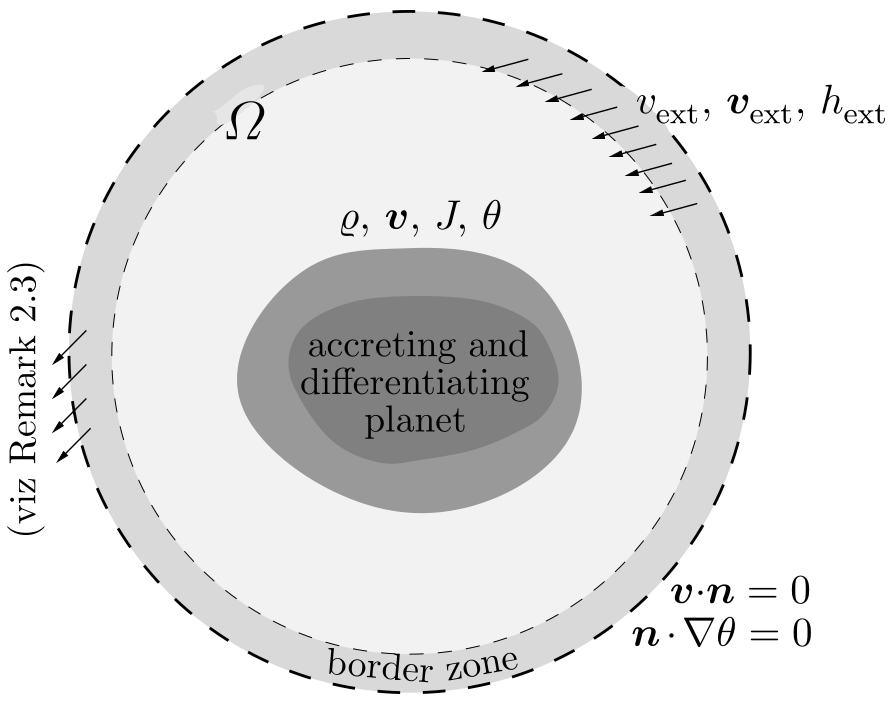}
%\end{my-picture}
\end{center}
\vspace*{-.6em}
\caption{
{\sl Schematic geometry of the fixed bounded domain $\varOmega$
and a border zone in the vicinity of its boundary $\varGamma$ where the
inflowing material from outer space is assumed to occur (and possibly
also to disappear), while the planetary accretion and differentiation
due to self-gravitation takes place  far from the border
zone.  }
}
\label{fig-accretion}
\end{figure}

 As regards the modelization of the incoming material, a natural
approach would be that of prescribing the material influx  from
outer space at the boundary $\varGamma=\partial\varOmega$. Although natural, this
setting would bring major analytical
difficulties, cf.\ \cite{FeiNov22MOFS},  which seem to preclude the
possibility of obtaining a comprehensive theory, even in the more regular case of 
higher-grade nonsimple viscous media, see in Remark~\ref{rem-hyper} below.
 We are hence to adopt a simplification, by fixing a-priori a specific bulk
 subdomain  of $\varOmega$  in the vicinity of $\varGamma$ (which we call
\ttt a \eee {\it border zone} in the following) far from the central region
where the planet is actually accreting, where we model the incoming material as
a given bulk source, see Figure~\ref{fig-accretion}. In particular, we
consider   such a bulk source to continuously provide material at some
given rate of influx of mass density, volume, velocity,  and  heat.

 A second specific aspect of the modelization is the description of the gas-to-fluid
thermo-mechanical transition at the planet surface. This corresponds to
a localized transition in terms of density or, equivalently, actual
volume. Specifically, we model such phase transition by prescribing
an abrupt variation (discontinuity) of the pressure in terms of the
density (or actual volume), see Figure \ref{PT-accretion}.

Albeit various computational simulations are reported in the 
geophysical literature,  see,  e.g., \cite{BadWal15EEAD,Cham04PAIS},
 a  specific  thermodynamically consistent  continuum-mechanical model taking into
account accretion and differentiation seems to be still missing. In
fact, to model such coupled effects was indicated as one of the open 
(numerical) challenges  in the field  in \cite[p.~417]{Gery19INGM}.
Let us mention that viscous effects are  %while viscosity is
standardly considered in dense planets (see, e.g.,
\cite{Gery19INGM,GerYue07RCMM,MulWei16SDES}) where viscosity is related in
particular with frictional (including tidal) heating.  In addition, the
very diluted extraplanetary dust may be also considered viscous. Indeed, such viscosity is at the basis of the  
so-called dynamical (Chandrasekhar) friction due to gravitational
self-interaction, cf., e.g., \cite{Cham04PAIS,IdaMak92NSGI,OBMoLe06TPFS}. This
justifies our setting where the whole material system in $\varOmega$
is assumed to be viscous. Notably, this assumption allows to obtain
some analytical results, as mentioned in Section~\ref{sec-anal}.

The plan of this paper is the following. In Section~\ref{sec-single},
the standard Eulerian mechanics and thermomechanics for viscoelastic
fluids is generalized to open systems in contact with the outer space
by prescribing a supply of mass (including its volume), momenta, and energies. 
In particular, the onset of fictitious forces related to the orbiting
of the system $\varOmega$ around a distant orbital center is
discussed. In Section~\ref{sect-two-comp}, this basic
single-component model is extended for a two-component variant
for materials with different mass densities and viscoelastic
properties. In Section~\ref{sec-anal}, we  discuss the stability of
the model in terms of possible a-priori estimates. By 
assuming a sufficiently fast blow-up of the stored energy,  we show that 
self-gravitational (non-relativistic) collapse at finite time  is
excluded. Eventually, the prospects of a mathematical analysis of a
multipolar variant  of the model are  briefly commented.

\section{A single-component model}\label{sec-single}
%        ~~~~~~~~~~~~~~~~~~~~~~~

 We start by discussing a single-component model, before moving to
the two component  metal/silicate model in
Section~\ref{sect-two-comp}.  We   anticipate the main
notation used in this paper, as in the following table:

\begin{center}
\fbox{
\begin{minipage}[t]{0.39\linewidth}\small\smallskip
  $\yy$ deformation (in m),\\
$\FF$ deformation gradient,\\
  $\vv$ velocity (in m/s),\\
$\varrho$ mass density (in kg/m$^3$),\\
  $\theta$ temperature (in K),\\
  $\TT$ conservative Cauchy stress (in Pa),\\
$\DD$ viscous (dissipative) stress (in Pa),\\
  $v_{\rm ext}^{}$ volume rate inflow (in 1/s),\\
  $\vv_{\rm ext}^{}$ velocity of inflow (in m/s),\\
  $h_{\rm ext}^{}$ heat rate inflow (in W/m$^3$),\\
  $c$ heat capacity (in %J/(kg\,K)
  Pa/K),\\
$\xi$ dissipation rate (in W/m$^3$=Pa/s),\\
$(\cdot)'$ (partial) derivative,
\smallskip \end{minipage}
\begin{minipage}[t]{21.5em}\small\smallskip
$\bm\xi=\yy^{-1}$ the return mapping (in m),\\
$J=\det\FF$ Jacobian,\\
$\psi$ free energy (in J/m$^3$=Pa),\\
$\varrho_0$ given referential mass density (in kg/m$^3$),\\
$\eta$ entropy (in Pa/K),\\
$\W$ heat internal energy (in Pa),\\ 
 $r_{\rm ext}^{}$ mass density rate inflow (in kg/(m$^3$s)),\\
$\ee(\vv)=\frac12\nabla\vv^\top\!\!+\frac12\nabla\vv$ small strain rate (in s$^{-1}$),\\
  $\GRAVPOT$ gravitational potential (in J/kg),\\
$\GRAVCONST$  gravitational constant  (in m$^3$kg$^{-1}$s$^{-2}$),\\
  $\bm g=-\GRAVCONST\nabla\GRAVPOT$ gravity acceleration (in m/s$^{2}$),\\
 $\kappa$ thermal conductivity (in W/(m\,K)),
 \\
  $(\cdot)\!\DT{^{}}$ convective time derivative.
  \smallskip \end{minipage}
}\end{center}

\vspace{-.9em}

\begin{center}
{\small\sl Table\,1.\ }
{\small
Summary of the basic notation;
$\GRAVCONST\doteq 6.674{\times}10^{-11}$m$^3$kg$^{-1}$s$^{-2}$.}
\end{center}

\subsection{Kinematics at large strains in brief}\label{sec-kinem}
%           ~~~~~~~~~~~~~~~~~~~~~~~~~~~~~~~~~~~

Let us start by recalling some basic notions from the general
theory of large deformations in continuum mechanics, limiting
ourselves to a minimal frame, to serve the sole
purpose of presenting the model, and referring the reader, e.g., to  
\cite{GuFrAn10MTC,Mart19PCM} for additional material.   

Assume to be given the {\it deformation} $\yy:I{\times}\varOmega\to\R^3$,
where $I=[0,T]$ and $T>0$ is some terminal time. For all given times
$t\in I$, the deformation maps the reference configuration
$\varOmega\subset\R^3$ of the deformable medium to its {\it actual} configuration
$\yy(t,\varOmega)$, a subset of the physical space $\R^3$. In what
follows, we indicate {\it referential coordinates} by $\XX\in\varOmega$
and {\it actual coordinates} by $\xx\in\R^d$.  We will work in the
Eulerian setting, i.e.,  in the  actual configuration. Unless otherwise
stated, all quantities are assumed to be actual and to be dependent on
the actual coordinates $\xx$, in case they are space dependent. All
other cases will be explicitly pointed out.

By assuming $\XX\mapsto\yy(t,\XX)$ to be globally invertible,
cf.\ Remark~\ref{rem-invert} below, we indicate its inverse by
$\bm{\xi}(t,\cdot)=\yy^{-1}(t,\cdot):\yy(t,\varOmega)\to\varOmega$.
Standardly, $\bm\xi$ is called {\it reference mapping} or {\it return mapping}.

The {\it referential velocity}  reads  $\frac{\d}{\d
  t}\yy( t, \XX)$  and   
the {\it  actual  velocity} is $\vv=\vv(\xx)=[\frac{\d}{\d t}\yy](\bm\xi(\xx))$.
The  actual  velocity $\vv$ is used to define the {\it material
derivative} $\DT q(t,\xx)$ of a scalar Eulerian quantity $q=q(t,\xx)$ as
$$
\DT q=\pdt q+\vv\Cdot\nabla q\,.
$$
For a vector-valued $\qq$, the material derivative is defined componentwise and
written as $\DT\qq=\pdt{}\qq+(\vv\Cdot\nabla)\qq$. With this notation, the return
mapping is evolving by the simple transport equation
\begin{align}\label{eq:xi0}
\DT{\bm\xi}=\bm0\,.
\end{align}

Given the deformation $  
\yy(t,\XX)$,  we define the
{\it deformation gradient} by
$\FF=\FF( t,\xx)=\nabla_{\XX}^{}\yy( t,  \bm\xi(\xx))$.
For brevity, we will write the composition with the return mapping $\bm\xi$ as an
upper index, i.e., for any  field $Q=Q(\XX)$ given in terms of the
referential coordinates,  the corresponding
field $Q(\bm\xi(\xx))$  in actual coordinates  will be 
denoted by  $Q^{\bm\xi}=Q^{\bm\xi}(\xx)$.

The deformation gradient evolves according the equation
$\DT\FF=(\nabla\yy)\FF$ and is related to the return mapping by the algebraic
relation
\begin{align}
  \FF=(\nabla\bm\xi)^{-1}\,.
\end{align}
In particular, the {\it Jacobian} $J:=\det\FF=1/\!\det(\nabla\bm\xi)$ evolves by the
equation
\begin{align}\label{J-evolv}
\DT J=({\rm div}\,\vv)J\,.
\end{align}

\subsection{Mechanics and thermomechanics in brief}\label{sec-thermodyn}
%           ~~~~~~~~~~~~~~~~~~~~~~~~~~~~~~~~~~~~~~

The well-established concept of {\it hyperelasticity} relies on 
the   existence of a stored energy  density  $\varphi=\varphi(\FF)$.
We consider this energy  density to be  expressed in physical
units J/m$^3$=Pa,  where the
volume is measured in the reference configuration. 
The corresponding {\it Cauchy stress} reads
\begin{equation}
\TT(\FF)=\frac{\varphi'(\FF)\FF^\top\!}{\det\FF}\,. \label{cau}
\end{equation}
An alternative option would be to prescribe the referential stored energy in
J/kg,  so that the Cauchy stress would take the form 
$\TT(\FF)=\varrho\varphi'(\FF)\FF^\top$,  involving 
the mass density $\varrho$,  see, e.g., 
\cite{MarHug83MFE} or \cite{GuFrAn10MTC} for both options,
respectively.  Here, we follow the former choice, leading to
\eqref{cau}, as the latter one seems less practical in presence of incoming
material. We consider $\varphi$ to be spatially
homogeneous. The heterogeneity of the system will then arise from the
distinct densities of the diluted material (gas and dust) and of the
compact planetary core. 

The mass density $\varrho$ standardly satisfies the continuity equation
\begin{align}\label{cont-eq}
\pdt\varrho+{\rm div}(\varrho\vv)=0\,.
\end{align}
Assuming the initial condition $\yy(0,\XX)=\XX$ so that
$\FF(0,\xx)=\bbI$  (identity tensor)  and 
considering also the initial condition $\varrho(0,\xx)=\varrho_0(\XX)$
 for some given initial density in the reference configuration $\XX\in\varOmega \mapsto \varrho_0(\XX)$, 
 relation  \eq{cont-eq} is equivalent to the algebraic equation
\begin{align}\label{cont-eq+}
\varrho=\frac{\varrho_0{\circ}\bm\xi}{\det\FF}\,.
\end{align}
This relation relies on \eq{J-evolv}. Conversely,
\eq{cont-eq} can be deduced from \eq{cont-eq+} when using \eq{J-evolv}.

 Planetary accretion and core-mantle formation occur on very long
time scales, namely, billions of years. On such time scales, solid-type
rheologies are inadequate  to describe material response  and  one   should 
rather  resort to  a fluidic rheology in deviatoric 
variables. We choose here the {\it Newton} {\it rheology}
(also called {\it Stokes} or, in the convective variant, the {\it
  Navier-Stokes}). Note that the Newton rheology does not predict the
propagation of shear waves, as opposed to Jeffreys' rheology, which
is often used in geophysical models.  We find nevertheless the choice
of the Newton rheology well justified, as shear weaves are negligible on
the time scale of planetary evolution. 

Fluids are characterized by  vanishing shear elastic resistance. The stored
energy then takes the form $\varphi(\FF)=\phi(J)$, cf., e.g.,
\cite[p.\,10]{MarHug83MFE}.  We assume $J$ to be positive, to
model local orientation preservation.  The conservative
Cauchy stress then reduces to  the   {\it pressure tensor}
\begin{align}\label{T-xi-J/kg-fluid}
\TT(\FF)=
\frac{\varphi'(\FF)\FF^\top\!}{ J }=\frac{\phi'( J )\det{}'(\FF)\FF^\top\!\!}{ J }
=\frac{\phi'( J )({\rm Cof}\,\FF)\FF^\top\!\!}{ J }=\phi'(J)\bbI\,
\end{align}
where the last equation  follows by Cramer's rule. From now on, we will
consider this fluidic ansatz, involving only $J$ instead of the whole tensor
$\FF$.

Besides $\vv$ and $\FF$, another important variable is  the 
temperature $\theta$.  We indicate {\it free energy} density  by 
$\psi=\psi(J,\theta)$,  to be  considered again in Pa, 
and the dissipation potential  density 
$\DIS=\DIS(J,\theta;\ee)$ in Pa/s, with $\ee$  being 
a placeholder for $\ee(\vv)$.  Let us specify that the densities
$\psi$ and $\DIS$ are assumed to be computed with respect to
undeformed, referential volumes.  Similarly as in \eq{T-xi-J/kg-fluid},
the Cauchy stress for the whole free energy  is $\psi_J'(J,\theta)\bbI$.
The momentum equilibrium equation then reads as
\begin{align}
\!\!\pdt{}(\varrho\vv)+{\rm div}(\varrho\vv{\otimes}\vv)
={\rm div}(\TT{+}\DD)+\ff\ &\text{ with }\,
\TT=\psi_J'(J,\theta)\bbI\
\text{ and }\,\DD=\DIS_{\ee}'(J,\theta;\ee(\vv)),\!\!
\label{momentum-eq}\end{align}
where $\ff$ is a bulk force,  to be later specified in Section
\ref{sec:gravity} and involving gravitational and
referential contributions.

 The {\it entropy} is written as
$\eta=-\psi_\theta'(J,\theta)/J$. The factor $J$ at the denominator reflects
the fact that, again, we measure volumes in the reference
configuration. The {\it entropy equation}  then reads 
\begin{align}\nonumber
 \pdt\eta+{\rm div}\big(\vv\eta\big)
=\frac{\xi-{\rm div}\,{\bm j}}\theta\ \ \ \ &\text{ with }\ 
  \jj=-\kappa(J,\theta)\nabla\theta \ \ \ \ \text{(i.e.,\ Fourier's law)}
\\&
\text{ and }\ \ \xi=\DD{:}\ee(\vv)=\DIS_{\ee}'(J,\theta;\ee(\vv)){:}\ee(\vv)\,,
\label{entropy-eq}\end{align}
where $\kappa=\kappa(J,\theta)$ is the {\it thermal conductivity coefficient} in W/(m$^2$K)
and $\xi=\xi(J,\theta,\ee(\vv))$ is the {\it heat production rate} due to mechanical
viscosity in Pa/s. Substituting  the expression for  $\eta$
into the entropy equation  gives 
\begin{align}\nonumber
\pdt\eta+{\rm div}\big(\vv\eta\big)=\DT\eta+\eta\,{\rm div}\,\vv
&=-\ \DT{\overline{\!\!\Big(\frac{\psi_\theta'(J,\theta)}J\Big)\!\!}}\ 
-\frac{\psi_\theta'(J,\theta)}J\,{\rm div}\,\vv
\\&
=-\frac{\psi_{\theta\theta}''(J,\theta)}J\,\DT\theta+\Big(
\frac{\psi_\theta'(J,\theta)}{J^2}{-}\frac{\psi_{J\theta}''(J,\theta)}J\Big)\,\DT J
-\frac{\psi_\theta'(J,\theta)}J\,{\rm div}\,\vv
\nonumber\\&\!\!\stackrel{\scriptsize\eq{J-evolv}}{=}\!
-\frac{\psi_{\theta\theta}''(J,\theta)\!}J\,\DT\theta-\psi_{J\theta}''(J,\theta)\,{\rm div}\,\vv\,.
\label{entropy-calculus}\end{align}
Thus, \eq{entropy-eq}  can be rewritten as  the {\it heat-transfer equation}
\begin{align}
c(J,\theta)\DT\theta-{\rm div}\big(\kappa(J,\theta)\nabla\theta\big)=
\xi+\theta\psi_{J\theta}''(J,\theta)\,{\rm div}\,\vv\quad
\text{ with }\ \ c(J,\theta)=-\theta\frac{\psi_{\theta\theta}''(J,\theta)}J\,,
\label{heat-eq}\end{align}
where $c=c(J,\theta)$ is the {\it heat capacity}.
The {\it internal energy} is given by the {\it Gibbs relation}
$\psi+\theta\eta$. We will use the (uniquely determined) split:
\begin{align}
\linesunder{\psi(J,\theta)-\theta\psi_\theta'(J,\theta)}{internal}{energy}\!\!
=\!\!\!\!\!\!\linesunder{\psi(J,0)}{$\ =:\phi(J)\;$stored}{\ \ \ \ \ \ \ \ \ \ \ energy}
\!\!\!\!\!\!\!\!+\!\!
\linesunder{\psi(J,\theta)-\theta\psi_\theta'(J,\theta)
-\psi(J,0)}
{thermal part of the}{referential internal energy}\!\!\!.
\label{split+}
\end{align}
 Working in the actual configuration, we will deal with the
quantities  
\begin{align}
\linesunder{\frac{\psi(J,\theta)-\theta\psi_\theta'(J,\theta)}J}{actual internal}{energy}\!\!\
=\!\!\!\!\!\!\!\!\linesunder{\frac{\psi(J,0)}J}{actual}{stored energy}\!\!\!\!\!\!\!\!
+\!\!\!
\linesunder{\frac{\psi(J,\theta)-\theta\psi_\theta'(J,\theta)-\psi(J,0)}J}
{$=:\W(J,\theta)\:$ thermal part of the}{actual internal energy}\!\!\!;
\label{split++}
\end{align}
where the term {\it actual} above is meant to point out that densities
are computed with respect to actual volumes. Note  that $\W(J,0)=0$. 

In terms of  $\W$ defined in \eq{split++}, the heat equation \eq{heat-eq}
can be written in the so-called {\it enthalpy formulation}
\begin{align}
&\hspace*{-0em}\pdt\W
+{\rm div}\big(\vv\W{-}\kappa(J,\theta)\nabla\theta\big)=\xi
+\big(\psi_{J}'(J,\theta){-}\phi'(J)\big)\,{\rm div}\,\vv\ 
\label{thermoviscoelastic-fluid-kg/J}\end{align}
with $\W$ from \eq{split++}. Here, using the shorter notation for the coupling
part of the free energy $\COUPLING(J,\theta):=\psi(J,\theta)-\phi(J)$, so that
$\W=(\COUPLING(J,\theta){-}\theta\COUPLING_\theta'(J,\theta))/J$, we relied on the
calculus
\begin{align}\nonumber
\pdt\W&+{\rm div}(\vv\W)=\DT\W+\W\,{\rm div}\,\vv
=\DT{\overline{\!\!\Big(\frac{\COUPLING(J,\theta){-}\theta\COUPLING_\theta'(J,\theta)}J\Big)\!\!}}
\ +\frac{\COUPLING(J,\theta){-}\theta\COUPLING_\theta'(J,\theta)}J\,{\rm div}\,\vv
\\&\nonumber\
=\bigg(\Big(\frac{\COUPLING(J,\theta)}J\Big)_J'{-}
\theta\Big(\frac{\COUPLING_\theta'(J,\theta)}J\Big)_J'\bigg)\DT J
-\theta\frac{\COUPLING_{\theta\theta}''(J,\theta)}J\DT\theta
+\frac{\COUPLING(J,\theta){-}\theta\COUPLING_\theta'(J,\theta)}J\,{\rm div}\,\vv
\\&\nonumber
\stackrel{\scriptsize\eq{J-evolv}}{=}\!J\bigg(\Big(\frac{\COUPLING(J,\theta)}J\Big)_J'{-}
\theta\Big(\frac{\COUPLING_\theta'(J,\theta)}J\Big)_J'\bigg)\,{\rm div}\,\vv
-\theta\frac{\COUPLING_{\theta\theta}''(J,\theta)}J\DT\theta
+\frac{\COUPLING(J,\theta){-}\theta\COUPLING_\theta'(J,\theta)}J\,{\rm div}\,\vv
\\&
\ =c(J,\theta)\DT\theta+\COUPLING_J'(J,\theta)\,{\rm div}\,\vv
-\theta\COUPLING_{J\theta}''(J,\theta)\,{\rm div}\,\vv
\label{enthalpy-calculations}\end{align}
with $c(J,\theta)$ from \eq{heat-eq}.
This shows that \eq{heat-eq} is indeed equivalent to
\eq{thermoviscoelastic-fluid-kg/J}.

\subsection{ Gravitational and fictitious forces}\label{sec:gravity}

\begin{figure}[ht]
%\begin{center}
%\begin{my-picture}{6}{6.5}{fig1}
\psfrag{o}{\small$\bm{\omega}$}
\psfrag{sun}{\footnotesize\begin{minipage}[t]{10em}a central\\[-.2em]star (Sun)\end{minipage}}
\psfrag{e1}{\footnotesize $\bm{e}_1$} 
\psfrag{e2}{\footnotesize $\bm{e}_2$} 
\psfrag{e3}{\footnotesize $\bm{e}_3$} 
\psfrag{hat-e1}{\footnotesize $\widehat{\bm{e}}_1$} 
\psfrag{hat-e2}{\footnotesize $\widehat{\bm{e}}_2$} 
\psfrag{hat-e3}{\footnotesize $\widehat{\bm{e}}_3$} 
\psfrag{W}{\large$\varOmega$}
\psfrag{G}{\large$\varGamma$}
\psfrag{viz Remark 2.2}{\footnotesize (viz Remark~2.3)}
\psfrag{border}{\footnotesize border}
\psfrag{zone}{\footnotesize zone}
\psfrag{V=0}{}\psfrag{v=0}{\small $\vv{\cdot}\nn=0$} 
\psfrag{v,J}{\small $\varrho,\,\vv,\,J,\,\theta$}
\psfrag{r-ext}{\small $v_{\rm ext}^{},\,\vv_{\rm ext}^{},\,h_{\rm ext}^{}$}
\psfrag{flux=0}{\small $\nn\Cdot\nabla\theta=0$}
\psfrag{planet}{\footnotesize\begin{minipage}[t]{17em}the orbiting accreting\\[-.2em]\hspace*{1.5em}planet surrounded by\\[-.2em]\hspace*{2.5em}a fictitious zone (Fig.\,1)\end{minipage}}
%\hspace*{1em}\includegraphics[width=29em]{tr-us-accretion3.eps}
\hspace*{3em}\includegraphics[width=36em]{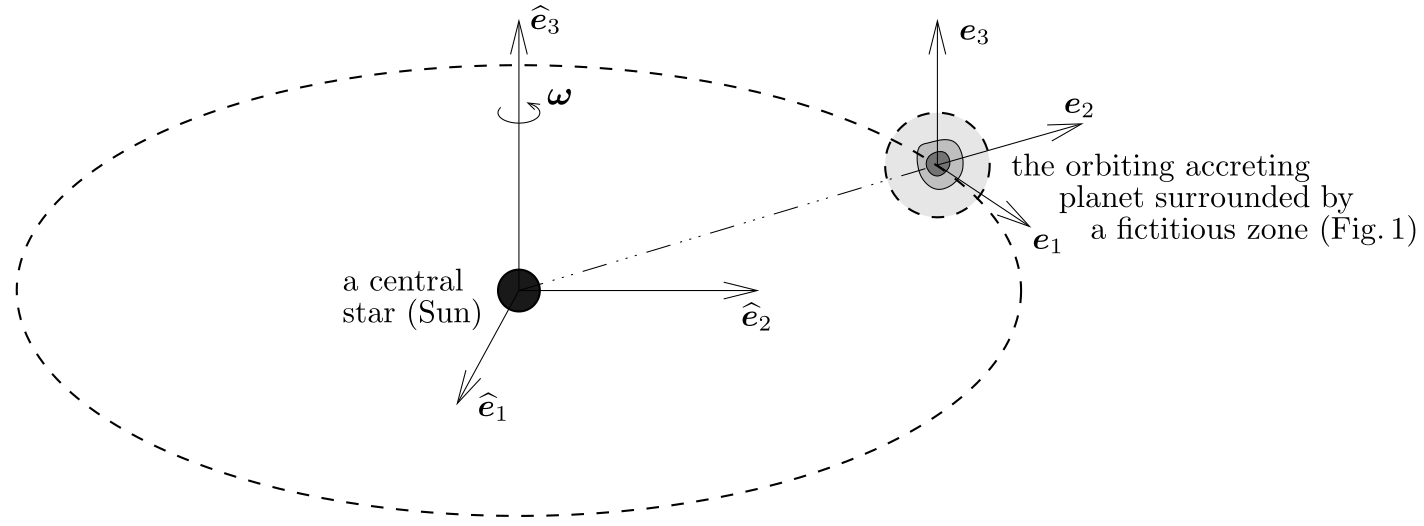}
%\end{my-picture}
%\end{center}
\vspace*{-.6em}
\caption{
  {\sl  The system from Figure~\ref{fig-accretion} is orbiting
    around a distant orbital center. Such motion is encoded by
    assuming that the coordinate frame  $(\bm{e}_1,\bm{e}_2,\bm{e}_3)$
    to be in uniform circular motion with respect to a fixed inertial
    referential frame
    $(\widehat{\bm{e}}_1,\widehat{\bm{e}}_2,\widehat{\bm{e}}_3)$  with
the angular velocity $\bm{\omega}$. Correspondingly, the referential
forces occur in the frame   $(\bm{e}_1,\bm{e}_2,\bm{e}_3)$.
}
}
\label{fig-orbiting}
\end{figure}

Let us now come to modeling gravitation. 
We assume that the system is orbiting around a very distant orbital
center, where a large mass is concentrated. An inertial reference
frame is attached to this large mass, so that $\xx$ corresponds to an
orbiting and rotating noninertial reference frame instead. More precisely,
we assume that the orbital motion of the system is regular circular with
axis ${\bf e}_3$ and
that the noninertial reference frame rotates along the orbit in such a
way that the axis $-{\bf e}_2$ is always directed toward the orbital
center,  see Figure \ref{fig-orbiting}. 

It is convenient to introduce two reference frames. 
Let $(\haz {\bm
  e}_1,\haz {\bm
  e}_2, \haz {\bm e}_3)$ indicate an inertial reference frame centered
in the orbital center and indicate the corresponding variable as $\haz
\xx=( \haz x_1, \haz x_2,  \haz x_3)$. At the same time, we consider
the 
noninertial reference frame $({\bm e}_1,{\bm e}_2, {\bm e}_3)$
rotating with respect to $(\haz {\bm
  e}_1,\haz {\bm e}_2,\haz {\bm e}_3)$ with axis
$\haz {\bm e}_3$, same center, and with constant scalar angular
velocity $\omega$. We indicate the
corresponding variable by $\xx=( x_1, x_2,  x_3)$. More precisely, we
consider the coordinate transformation $\haz \xx = R_t(\xx)$ by letting 
$$ \haz x_3 =  x_3,\quad \haz x_1= x_1\cos(\theta(t)) - x_2
\sin(\theta(t)), \quad \haz x_2 = x_1\sin(\theta(t)) +
x_2\cos(\theta(t))$$
with $\theta(t)= \omega t  $ for $t\in \R$. In particular, one has that
${\bm e}_3=\haz {\bm e}_3$, so that the angular velocity vector
$\oomega=(0,0,\omega)$ has the same expression in both reference
frames.

The density of the system in inertial coordinates reads $\haz
\varrho(t,\haz \xx) = \varrho(t,R^{-1}_t(\haz \xx)) = \varrho(t,\xx)$. We
assume that the orbital center $0$ is the barycenter of a large mass $\haz
\varrho_*(t,\haz \xx)$. In case of the
Sun-Earth system, $\haz \varrho_*$ is the mass distribution of the Sun
and $\haz \varrho$ is that of the Earth. The two mass distributions
$\haz \varrho_*$ and
$\haz \varrho$ are assumed to be very distant and to be  both
compactly supported, uniformly with respect to time. 

The gravity acceleration
$\haz\gg(t,\haz\xx)=-\nabla_{\haz \xx}\haz\GRAVPOT(t,\haz\xx)$ is defined via the
gravitational constant $\GRAVCONST$ from the {\it gravitational potential}
$\haz\GRAVPOT$ governed for all times by the Poisson equation
\begin{align}\label{Poisson}
\Delta_{\haz \xx}\haz\GRAVPOT=4\uppi\GRAVCONST(\haz \varrho_* +\haz \varrho) \ \ \text{ on }\ \R^3\,.
\end{align}
In the physically relevant three-dimensional situation, the unique solution of
\eq{Poisson} satisfying the condition
$\lim_{|\haz \xx|\to\infty}\haz V(t,\haz \xx)=0$ can be written in
terms of the Green function of the Laplace operator as
\begin{align}\label{V-integral}
\haz \GRAVPOT(t,\haz \xx)=-\GRAVCONST\!\int_{\R^3}\!\frac{\haz
  \varrho_*(t,\widetilde\xx)+\haz \varrho(t,\widetilde\xx)}{|\haz \xx{-}\widetilde\xx|}\,\d\widetilde\xx\,.
\end{align}
Note that the potential $\haz \GRAVPOT$ is negative, following  the
geophysical convention (as opposed to the convention in mathematics, cf., e.g.,\
\cite{DFPS04GTWS,DuNeVa10GMCB}). Then we have
\ttt the mere gravitational acceleration \eee
\begin{align}\label{nabla-V-integral}
\haz {\bm g}(t,\haz \xx)=-\nabla_{\haz \xx}\haz \GRAVPOT(t,\haz \xx)
=\GRAVCONST\!\int_{\R^3}\!\frac{\big(\haz
  \varrho_*(t,\widetilde\xx)+\haz \varrho(t,\widetilde\xx)
  \big)(\widetilde\xx{-}\haz \xx)}{|\widetilde\xx{-}\haz \xx|^3}\,\d\widetilde\xx\,
\end{align}
so that the gravitational force density at $(t,\haz \xx)$ reads $\haz
\varrho(t,\haz \xx) {\bm g}(t,\haz \xx)$.

The \ttt overall \eee gravity acceleration relative to the uniformly rotating
noninertial  frame $({\bm e}_1,{\bm e}_2,{\bm e}_3)$ takes the form
$$
\haz{\bm g}-2\oomega{\times}\vv -\oomega{\times}(\oomega{\times}\haz\xx)\,,
$$
where $- 2  \oomega \times \vv $ is the Coriolis
acceleration (and $\vv$ denotes the velocity relative to the rotating
frame) and $-\oomega{\times}(\oomega{\times}\haz \xx )$ is the centrifugal
acceleration. Note that no Euler force $-\dot \oomega{\times}\haz \xx$ occurs,
as $\oomega$ is constant.

Our aim is now to reduce the model to the region $\varOmega$ and to the
noninertial coordinates $\xx$. To this end, we introduce some simplifying
assumption, which apply to the Sun-Earth system, in particular. To start with,
we assume that the mass of $\haz\varrho_*$ is much larger than that of
$\haz \varrho$. This allows to neglect the influence of $\haz\varrho$ on
$\haz \varrho_*$, so that we reduce ourselves to model the influence of the
latter on the former.

Secondly, we assume that the centrifugal forces are balanced by the
gravitational forces generated by the large mass distribution $\haz
\varrho_*$, so that both contributions can be omitted. This assumption
is reasonable in case the supports of the two mass distributions is much smaller
than their distance. In fact, in case the support of $\haz \varrho_*$ is 
small and independent of time one can approximate the gravity acceleration
generated by $\haz \varrho_*$ by
$ -{\GRAVCONST M_* \haz \xx}/{|\haz \xx|^3}$ where $M_*$ indicates the
mass $ 
M_*:=\int_{\R^3}\haz \varrho_*(t,\tilde \xx)\, \d \tilde \xx$.
On the other hand, if the support of $\haz \varrho$ is small, one can
approximate the centrifugal acceleration by
$-\oomega{\times}(\oomega{\times}\haz \xx)=\omega^2(\haz x_1,\haz x_2,0)$ by
$\omega^2 \haz\xx$, i.e., by modifying the ${\bf e}_3$ component. Let now
$D>0$ be the distance between the two
supports. If $D$ is large with respect to the diameter of the support
of $\haz \varrho$, in order to quantify the acceleration on $\haz \varrho$ one can systematically replace
$\haz \xx $ by $D\haz \xx/|\haz \xx|$ in the expression of the
accelerations. In particular, the centrifugal and gravity
accelerations on $\haz \varrho$ are balanced if   
$$
{\bm 0}=\frac{GM_*\haz \xx}{D^3}\sss-\eee\frac{\omega^2D \haz \xx}{D}
$$
which holds for $\omega^2=GM_*/D^3$ (note that the sign of $\omega$ is
uninfluential). Let us mention that the above approximations are well
justified for the Sun-Earth system. In this case, the computed revolution period $2\pi/\omega$ corresponds
to 365.26 days, to be compared with 365.24 days duration of one year,
with a 0.005\% error.

In conclusion, under the above assumptions, by letting 
$ {\bm g} ( t,\xx) = \haz {\bm g} (t,R_t(\haz \xx))  $ the momentum equilibrium
equation \eqref{momentum-eq} reads
\begin{equation}
\pdt{}(\varrho\vv)+{\rm div}\big(\varrho\vv{\otimes}\vv
-\DIS_{\ee}'(J,\theta;\ee(\vv))\big)
+\nabla p= \varrho {\bm g}-2\varrho\oomega{\times}\vv\,.
\label{momentum-eq2}
\end{equation}

\subsection{Open thermodynamics of self-gravitating accretion}\label{sec-thermodyn-open} 
%           ~~~~~~~~~~~~~~~~~~~~~~~~~~~~~~~~~~~~~~~~~~~~~~~~~~

 We now proceed to extend the setting from
Section~\ref{sec-thermodyn} to allow for incoming material. It will
be instrumental to assume that the region $\varOmega$ is invariant
under deformation, i.e., $\yy(t,\varOmega)= \varOmega$ for all
times. This in particular implies that the normal component of the
velocity vanishes at the boundary, namely,  $\vv{\cdot} \nn= 0 $ on $\varGamma:=\partial \varOmega$. Such a no-flux boundary conditions does
not allow to model material influx across the boundary $\varOmega$. We
hence resort in modeling incoming material by identifying some
border region of $\varOmega$, where material is added, see
Figure~\ref{fig-accretion}. We believe this somewhat nonphysical
setting to be well motivated by the fact that the relevant portion of
$\varOmega$, where the system shows high density, may be assumed to be
far from this border region. For simplicity, let us ignore a possible
outflow of material, but see Remark~\ref{rem-outflowing},

In the border region, we assume that the incoming material has a given mass
density rate $r_{\rm ext}^{}=r_{\rm ext}^{}(t,\xx)$, a velocity
$\vv_{\rm ext}^{}=\vv_{\rm ext}^{}(t,\xx)$, and  a heat rate (power)
$h_{\rm ext}^{}=h_{\rm ext}^{}(t,\xx)$. These quantities may depend on time
and space.  On the other hand,  to avoid other technicalities, the incoming
material itself is assumed to feature time-and-space independent
thermo-viscoelastic properties. In particular, the referential mass density
$\varrho_0$ of the incoming material, its free energy $\psi$, and its viscosity
described by the dissipation $\DIS$ are assumed to not depend on time
nor  on space. In particular, the  mass conservation equation \eqref{cont-eq}
is extended as
\begin{equation}\label{eq:cont-eq}
\frac{\partial\varrho}{\partial t} + {\rm div} (\varrho \vv) = r_{\rm ext}\,.
\end{equation}

To prevent infinite compression, i.e., $J=0$, one can follow the approach in
\cite{MiRoSt23MGDC,Roub24TVSE}, where a positive lower
bound on $J$ is directly obtained from the corresponding evolution equation.
These however apply to the case of suitably regular velocities $\vv$, and
do not seem to be amenable in the current setting. We hence follow another
path and assume that the energy blows up for $J \to 0+$. At the same time,
we aim at conserving the algebraic relation
\eq{cont-eq+},  which in turn ensures the possibility of
controlling $\varrho$. Yet,  assuming  \eq{cont-eq+} and \eqref{eq:cont-eq}
requires  extending also the kinematic equation \eq{J-evolv} suitably. To
this goal, we assume that an incoming {\it volume rate} $-v_{\rm ext}^{}$  can
be associated with the incoming material, suitably related with 
$r_{\rm ext}^{}$,  cf.\ \eq{r-vs-v} below.
This  modifies  \eq{J-evolv} as
\begin{align}\label{J-evolv-ext}
\DT J=({\rm div}\,\vv-v_{\rm ext}^{})J\,.
\end{align}
 An  extension of a kinematic equation (in terms of the
whole Cauchy-Green tensor $\FF\FF^\top$)  of this kind  is used in
\cite{GaKoTr22VCHM} for a tumor-growth problem.

Altogether, this incoming bulk material  induces additional
contributions in the  right-hand-side 
 of  \eq{cont-eq} by $r_{\rm ext}^{\nothing}$  as in
\eqref{eq:cont-eq},  of  the momentum equation   \eq{momentum-eq2} by the incoming
momentum rate $r_{\rm ext}^{\nothing}\vv_{\rm ext}^{\nothing}$, 
of  the
evolution  equation for $J$ 
\eq{J-evolv}  as in \eqref{J-evolv-ext}, and  of  the entropy
equation  \eq{entropy-eq} by the incoming entropy rate
$v_{\rm ext}^{\nothing}\eta+h_{\rm ext}^{\nothing}/\theta$, i.e., 
\begin{align}
 \pdt\eta+{\rm div}\big(\vv\eta\big)
=\frac{\xi-{\rm div}\,{\bm j}+h_{\rm ext}^{\nothing}}\theta+v_{\rm ext}^{\nothing}\eta
\label{entropy-eq+}\end{align}
with $\jj$ and $\xi$ again from \eq{entropy-eq}. The extended kinematic
equation \eq{J-evolv-ext} induces a modification in the computation for
\eqref{entropy-calculus} leading to 
\begin{align}\nonumber
\pdt\eta+{\rm div}\big(\vv\eta\big)
&=-\frac{\psi_{\theta\theta}''(J,\theta)}J\,\DT\theta
+\Big(\frac{\psi_{\theta}'(J,\theta)}{J^2}-\frac{\psi_{J\theta}''(J,\theta)}J\Big)\,\DT J
-\frac{\psi_\theta'(J,\theta)}J\,{\rm div}\,\vv
\nonumber\\&\!\!\!\stackrel{\scriptsize\eq{J-evolv-ext}}{=}\!
-\frac{\psi_{\theta\theta}''(J,\theta)}J\,\DT\theta
+\Big(\frac{\psi_{\theta}'(J,\theta)}{J}-\psi_{J\theta}''(J,\theta)\Big)\,
\big({\rm div}\,\vv{-}v_{\rm ext}^{\nothing}\big)
-\frac{\psi_\theta'(J,\theta)}J\,{\rm div}\,\vv
\nonumber\\&=
-\frac{\psi_{\theta\theta}''(J,\theta)}J\,\DT\theta-\psi_{J\theta}''(J,\theta)\,{\rm div}\,\vv
+v_{\rm ext}^{\nothing}\Big(
\psi_{J\theta}''(J,\theta){-}\frac{\psi_{\theta}'(J,\theta)}{J}\Big)
\,.
\label{entropy-calculus+}\end{align}
The extended entropy equation \eq{entropy-eq+} then gives the heat-transfer
equation 
\begin{align}
\!c(J,\theta)\DT\theta-{\rm div}\big(\kappa(J,\theta)\nabla\theta\big)=
\xi+\theta\psi_{J\theta}''(J,\theta)\,{\rm div}\,\vv
%+v_{\rm ext}^{\nothing}\theta\Big(\frac{\psi_{\theta}'(J,\theta)}{J}{-}\psi_{J\theta}''(J,\theta)\Big)
\sss-v_{\rm ext}^{\nothing}\theta\psi_{J\theta}''(J,\theta)\eee
+h_{\rm ext}^{\nothing}.
\label{heat-eq+}\end{align}
\ttt{}with the heat capacity $c$ again from \eq{heat-eq}. \eee The modification of
the kinematic  equation \eq{J-evolv-ext} influences also the calculations for 
\eq{enthalpy-calculations} as
\begin{align}\nonumber
\pdt\W&+{\rm div}(\vv\W)=\DT\W+\W\,{\rm div}\,\vv
=\DT{\overline{\!\!\Big(\frac{\COUPLING(J,\theta){-}\theta\COUPLING_\theta'(J,\theta)}J\Big)\!\!}}
\ +\frac{\COUPLING(J,\theta){-}\theta\COUPLING_\theta'(J,\theta)}J\,{\rm div}\,\vv
\\&\nonumber\
=\bigg(\Big(\frac{\COUPLING(J,\theta)}J\Big)_J'{-}
\theta\Big(\frac{\COUPLING_\theta'(J,\theta)}J\Big)_J'\bigg)\DT J
-\theta\frac{\COUPLING_{\theta\theta}''(J,\theta)}J\DT\theta
+\frac{\COUPLING(J,\theta){-}\theta\COUPLING_\theta'(J,\theta)}J\,{\rm div}\,\vv
\\&\nonumber
\stackrel{\scriptsize\eq{J-evolv-ext}}{=}\!
J\Big(\frac{\COUPLING(J,\theta){-}\theta\COUPLING_\theta'(J,\theta)}J\Big)_J'
\big({\rm div}\,\vv{-}v_{\rm ext}^{\nothing}\big)
-\theta\frac{\COUPLING_{\theta\theta}''(J,\theta)}J\DT\theta
+\frac{\COUPLING(J,\theta){-}\theta\COUPLING_\theta'(J,\theta)\!}J\,{\rm div}\,\vv
\\&\nonumber
\ =
c(J,\theta)\DT\theta+\COUPLING_J'(J,\theta)\,{\rm div}\,\vv
-\theta\COUPLING_{J\theta}''(J,\theta)\,{\rm div}\,\vv
-J\Big(\frac{\COUPLING(J,\theta){-}\theta\COUPLING_\theta'(J,\theta)}J\Big)_J'
v_{\rm ext}^{\nothing}
\\&
\ =
c(J,\theta)\DT\theta+\COUPLING_J'(J,\theta)\,{\rm div}\,\vv
-\theta\COUPLING_{J\theta}''(J,\theta)\,{\rm div}\,\vv
+\big(\W
-\COUPLING_J'(J,\theta)+\theta\COUPLING_{J\theta}''(J,\theta)\big)v_{\rm ext}^{\nothing}\,,
\label{enthalpy-calculations-ext}\end{align}
where again $\W=(\COUPLING(J,\theta)-\theta\COUPLING_\theta'(J,\theta))/J$
with $\COUPLING(J,\theta)=\psi(J,\theta)-\phi(J)$ 
as in \eq{split++}.

Therefore, the enthalpy formulation of \eq{heat-eq+} reads as 
\begin{align}
\pdt\W+{\rm div}\big(\vv\W{-}\kappa(J,\theta)\nabla\theta\big)=\xi+
\hspace{-1em}
\lineunder{\big(\psi_{J}'(J,\theta){-}\phi'(J)\big))
\big({\rm div}\,\vv{-}v_{\rm ext}^{\nothing}\big)}{$=\COUPLING_J'(J,\theta)\DT J/J$}\hspace{-1em}
+
%\Big(\hspace{-.7em}\lineunder{\!\W\,{+}\,\theta\frac{\psi_{\theta}'(J,\theta)\!}J}{$=\COUPLING(J,\theta)/J$}\hspace{-.7em}\Big)v_{\rm ext}^{\nothing}
\ttt\W v_{\rm ext}^{\nothing}\eee+h_{\rm ext}^{\nothing}\,.
%\nonumber\\[-1.7em]
\label{thermoviscoelastic-fluid-kg/J+}\end{align}

\begin{figure}[ht]
\begin{center}
\psfrag{JT}{\scriptsize $1$ 
}
\psfrag{pT}{\scriptsize $p_0$}
\psfrag{p}{\scriptsize $p$}
\psfrag{J}{\scriptsize $J=\varrho_0/\varrho$}
\psfrag{J0}{\scriptsize $J_0$}
\psfrag{f}{\scriptsize $\phi$}
\psfrag{f,}{\scriptsize $p=-\phi'(J)$}
\psfrag{b}{\scriptsize $b$}
\psfrag{compression}{\scriptsize compression}
\psfrag{1}{\scriptsize 1}
\psfrag{0}{\scriptsize 0}
\psfrag{r0}{\scriptsize $\varrho_0^{}$}
\psfrag{blow-up}{\scriptsize\begin{minipage}[t]{11em}blow-up faster\\[-.2em]than $1/J$ for $J\to0+$\end{minipage}}
\psfrag{growth}{\scriptsize\begin{minipage}[t]{11em}growing slower\\[-.2em]than $\sqrt{p}$ for $p\to\infty$\end{minipage}}
\psfrag{straight segment}{\scriptsize\begin{minipage}[t]{11em}straight\\[-.3em]\hspace*{-.2em}segment\end{minipage}}
\psfrag{planetary}{\scriptsize\begin{minipage}[t]{11em}planetary\\[-.2em]\hspace*{0em}compression\\[-.2em]\end{minipage}}
\psfrag{disc}{\scriptsize\begin{minipage}[t]{11em}accretion disc\\[-.2em]\hspace*{0em}around the planet\\[-.2em]\end{minipage}}
\psfrag{transition}{\scriptsize\begin{minipage}[t]{11em}transition due\\[-.2em]\hspace*{1em}to impacts\\[-.2em]\end{minipage}}
\includegraphics[width=34em]{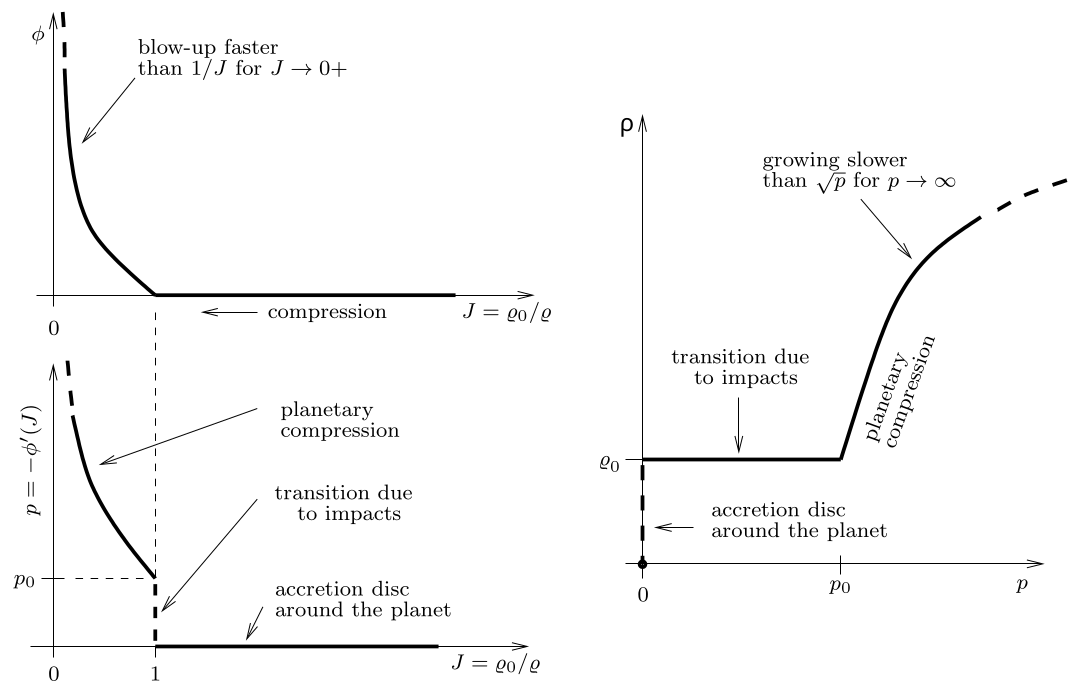}
\end{center}
\vspace*{-.6em}
\caption{
{\sl Schematic illustration of  the  transition from  the diluted
material (gas and dust) to the high-density viscoelastic fluid material of
the planetary core. 
This transition is  modeled by the kink in the graph of the 
bulk stored energy $\phi=\phi(J)$. 
}
}
\label{PT-accretion}
\end{figure}

An important point is the description of the thermo-mechanical
transition from the low-density dilute phase of the medium in areas
far  from  the planet to a high-density viscoelastic fluid after the impact
on the surface of the accreting planet.
This surface is  of course evolving  with time. Note that we are
not requiring this surface to be connected. In principle, the
accretion process might lead to the development of a binary 
planetary system  whose
two parts are orbiting  with respect to  each other.  This
is for instance the case the Pluto-Charon or the
Orcus-Vanth system. Moreover, our setting covers the case of a jovian-type
planet with rings. We assume that the mass density  at  the 
surface of the accreting planet, which essentially determines
the transition zone where impacts occur,
is considered  to be the  referential density $\varrho_0$ as in 
Figure~\ref{PT-accretion}.  The   sudden decrease  in  volume during such
impact-like transition leads
to adiabatic heating.  This,  together with  an abrupt  velocity variation
causing dissipation of mechanical energy through viscosity, leads
to  an  increase of temperature at the surface.

The self-gravitating open thermodynamical system is composed of
\eqref{eq:cont-eq}, \eq{momentum-eq2} expanded
by the momentum rate $r_{\rm ext}^{\nothing}\vv_{\rm ext}^{\nothing}$,
\eq{entropy-eq} expanded by $v_{\rm ext}^{\nothing}\eta+h_{\rm ext}^{\nothing}/\theta$
 and  resulting  into  \eq{thermoviscoelastic-fluid-kg/J+}, and of
\eq{J-evolv} expanded by $-v_{\rm ext}^{}J$  and  resulting  into
\eq{J-evolv-ext}. Altogether, it yields
the {\it extended} compressible Navier-Stokes-Fourier-Poisson system: 
\begin{subequations}\label{Euler-fluid-accret}\begin{align}
\label{Euler0-fluid-accret}&
\pdt\varrho+{\rm div}\big(\varrho\vv)=r_{\rm ext}^{\nothing}\,,\\
\nonumber
&\pdt{}(\varrho\vv)+{\rm div}\big(\varrho\vv{\otimes}\vv
-\DIS_{\ee}'(J,\theta;\ee(\vv))\big)
+\nabla p=r_{\rm ext}^{\nothing}\vv_{\rm ext}^{\nothing}
     +\varrho{\bm g}  - 2 \varrho \oomega{\times}\vv 
  \\[-.5em]&
  \hspace{12em}\text{with }\ p=-\psi_J'(J,\theta)\,,\ \ \ 
{\bm g}=\GRAVCONST\!\int_\varOmega\!\frac{\varrho(\widetilde\xx)(\widetilde\xx-\cdot)\!}{|\widetilde\xx-\cdot|^3}\,\d\widetilde\xx\,,
\label{Euler1-fluid-accret}
\\\nonumber
&\pdt\W+
{\rm div}\big(\W\vv-\kappa(J,\theta)\nabla\theta\big)
=\DIS_{\ee}'(J,\theta;\ee(\vv)){:}\ee(\vv)+
\big(\psi_{J}'(J,\theta){-}\phi'(J)\big)
\,\big({\rm div}\,\vv{-}v_{\rm ext}^{\nothing}\big)
\\&
\hspace{5.5em}
%+\DELETE{{\frac{\psi(J,\theta){-}\phi(J)\!}J\,v_{\rm ext}^{\nothing}}}
                                                \sss + w\,v_{\rm
                                                ext}^{\nothing} \eee +h_{\rm ext}^{\nothing}
\ \ \text{ with }\ \
\label{Euler4-fluid-accret}
\W=\frac{\psi(J,\theta)-\theta\psi_\theta'(J,\theta)-\phi(J)}J\,,
\\[-.3em]\label{Euler3-fluid-selfgravit-accret}
&\pdt J=({\rm div}\,\vv-v_{\rm ext}^{})J-\vv\Cdot\nabla J\,.
\end{align}\end{subequations}

 In Section~\ref{sec-anal} below, we investigate the {\it
  stability} of this system, namely, we prove that gravitational self-collapse
is excluded. To this aim, it is instrumental that $\varrho$  depends 
on $J$ through \eq{cont-eq+}. To this  end,  we pose the following
relation between incoming mass and volume rates
\begin{align}\label{r-vs-v}
r_{\rm ext}^{}=\frac{\varrho_0v _{\rm ext}^{}}{J}\,.
\end{align}
Let us note that \eq{Euler3-fluid-selfgravit-accret}, written as
\eq{J-evolv-ext}, is equivalent to
$\DT{\overline{(1/J)}}=(v_{\rm ext}^{}-{\rm div}\,\vv)(1/J)$. 
Hence,   by multiplying  by  $\varrho_0$ and  letting
  $\varrho=\varrho_0/J$ as in
\eq{cont-eq+}, we obtain $\DT\varrho=(v_{\rm ext}^{}-{\rm div}\,\vv)\varrho$.
This is  nothing but  \eq{Euler0-fluid-accret}, provided  that  \eq{r-vs-v} with $\varrho_0=\varrho J$
is taken into account.

Conversely, relation $\varrho=\varrho_0/J$ can be obtained from
\eq{Euler0-fluid-accret}  and  \eq{Euler3-fluid-selfgravit-accret}. Indeed,
assuming the existence of a function $f$ such that $\varrho=f(J)$, we have
\begin{align}
\frac{\DT\varrho}{\varrho}=\frac{f'(J)\DT J}{f(J)}=\frac{f'(J)J}{f(J)}\frac{\DT J}J
\stackrel{\scriptsize\eq{Euler3-fluid-selfgravit-accret}}{=}
\frac{f'(J)J}{f(J)}\big({\rm div}\,\vv-v_{\rm ext}^{}\big)
\end{align}
We  then recover \eq{Euler0-fluid-accret} if
$v_{\rm ext}^{}=r_{\rm ext}^{}/f(J)$ and $f'(J)J/f(J)=-1$, i.e., if $f$ solves the
ordinary differential equation $f'(J)=-f(J)/J$. Fixing $f( 1 )=\varrho_0$,
 the  equation  admits the  unique solution
$f(J)=\varrho_0/J$,  i.e., 
\eq{cont-eq+}.  In this case  \eq{r-vs-v} is
$r_{\rm ext}^{}=v_{\rm ext}^{}\varrho$.

Let us mention that an alternative  way to preserve \eq{cont-eq+}
would be that of modifying the  referential mass density $\varrho_0$, cf.\ 
\cite{CiAmMa12MTMP,CiaMau11EFSG,Epst15MCIR,EpsMau00TVGU}. Here, we
 do not to follow this alternative approach, as we prefer to stay
with the interpretation of $\varrho_0$ as initial condition.

\begin{remark}[{\sl The incoming volume rate}]\upshape
In order to offer a  possible interpretation of the occurrence of the incoming
volume rate  $v_{\rm ext}^{}$ let us recall the standard Kr\"oner-Lee-Liu
\cite{Kron60AKVE,LeeLiu67FSEP} {\it multiplicative
 decomposition} of $\FF=\FF_{\rm e}^{}\FF_{\rm s}^{}$.  Here,
 $\FF_{\rm e}^{}$  and  $\FF_{\rm s}^{}$ are classically interpreted
 as the elastic and the
 inelastic distortion, respectively. Depending on the context, the
 tensor $\FF_{\rm s}^{}$  is interpreted as a plastic or a creep distortion.
 In some models, it may
 correspond to some isotropic-swelling  distortion  and take the form 
 $\FF_{\rm s}^{}=(1{+}s)^{1/3}\bbI$, with $s$ denoting volume variation due
 to swelling or morphogenic growth, cf.\ e.g.\
 \cite{CheAna11TMCT,DuSoFi10TSMF} or \cite[Chap.12]{Gori17MMBG}. Then, letting
 $J_{\rm e}^{}=\det\FF_{\rm e}^{}$, we have $J=\det\FF=J_{\rm e}^{}(1{+}s)$. From 
$\DT J=(1{+}s)\DT J_{\rm e}^{}+J_{\rm e}^{}\DT s$ and \eq{J-evolv}, we have
  \begin{align*}
  \DT J_{\rm e}^{}=\Big({\rm div}\,\vv\,{-}\,\frac{\DT s}{1{+}s}\Big)J_{\rm e}^{}\,.
  \end{align*}
  This  is reminiscent of  \eq{J-evolv-ext}, with $J$ corresponding to 
  the elastic-distortion determinant  $J_{\rm e}$ and the incoming
  volume rate $v_{\rm ext}^{}$ to  the {\it volume-deviation rate}
  $\DT s/(1{+}s)$. Given $v_{\rm ext}^{}$, the volume deviation itself could
  then be reconstructed from the equation $\DT s=v_{\rm ext}^{}(1{+}s)$.
\end{remark}

\subsection{Mass, momentum, energy,  and entropy  balance behind \eq{Euler-fluid-accret}}
%           ~~~~~~~~~~~~~~~~~~~~~~~~~~~~~~~~~~~~~~~~~~~~~~~~~~~~~~~~~~~~~~~~~

To establish the balance of mass, momentum, energy,  and entropy  we must
prescribe boundary conditions. We consider the  domain $\varOmega$ with
boundary  $\varGamma:=\partial \varOmega$ to be
thermally insulated, impermeable, and offering no friction to motion 
(i.e.,  we assume  homogeneous Navier condition).
 In particular, we prescribe
 \begin{align}\label{Euler-fluid-accret-BC}
\vv{\cdot}\nn= 0 \,,\ \ \ \ \ \ \
\big[\DIS_{\ee}'(J,\theta;\ee(\vv))\nn\big]_\text{\sc t}=\bm0\,,
\ \ \ \text{ and }\ \ \
\nabla\theta{\cdot}\nn=0\ \ \ \ \text{ on } \ \varGamma\,.
\end{align}
where $[\cdot]_\text{\sc t}$ denotes the tangential component of a 
vector. The  Navier condition $\vv{\cdot}\nn =0$  ensures that the actual
deforming domain has the same shape as the referential one,  namely,
$\yy(t,\varOmega)=\varOmega$. 

The mass balance contains  a source term  due to the incoming mass
with the rate $r_{\rm ext}^{}$. By
integrating \eq{Euler0-fluid-accret} over $\varOmega$, we obtain
\begin{align}
\frac{\d}{\d t}\hspace{-.7em}\linesunder{\int_\varOmega\varrho\,\d\xx}
{total mass in}{the domain $\varOmega$}\hspace{-.7em}
=\!\!\linesunder{\int_\varOmega r_{\rm ext}^{\nothing}\,\d\xx}{incoming}{mass rate}\,.
\end{align}

The momentum balance  features a source term  due to the momentum of the
incoming material with the rate $r_{\rm ext}^{}\vv_{\rm ext}^{}$, and can be
obtained by integrating \eq{Euler1-fluid-accret} over $\varOmega$.
As for the self-gravitational acceleration, we can use the integral
formula \eq{nabla-V-integral}, which gives by symmetry
\begin{align}
\int_\varOmega\varrho\,{\bm g}\,\d\xx
=\GRAVCONST\!\int_\varOmega\!\int_\varOmega\frac{\varrho(\xx)\varrho(\widetilde\xx)
(\widetilde\xx{-}\xx)}{|\widetilde\xx{-}\xx|^3}\,\d\xx\d\widetilde\xx=\bm0\,.
\label{self-graviating-momentum}\end{align}
Altogether, we obtain the momentum balance:
\begin{align}
\frac{\d}{\d t}\hspace*{-.7em}\linesunder{\int_\varOmega\varrho\vv\,\d\xx}{total}{momentum}\hspace*{-.7em}
=\hspace*{-.7em}\linesunder{\int_\varOmega
r_{\rm ext}^{\nothing}\vv_{\rm
  ext}^{\nothing}\,\d\xx}{incoming}{momentum rate} \hspace{-.7em}-
  \hspace{-1em}\linesunder{\int_\varOmega2\varrho \oomega{\times}\vv\,\d\xx}
  {momentum rate from}{the Coriolis force}\hspace*{-.7em}.
\end{align}
 Note the occurrence of a momentum rate related to Coriolis
forces, which reveals the {\it external} nature of such fictitious
force. 

To proceed further, we rewrite the momentum equation \eq{Euler1-fluid-accret}
by using 
\begin{align}\nonumber
&\pdt{}(\varrho \vv)+{\rm div}(\varrho \vv{\otimes}\vv)
 =\varrho\pdt \vv+\pdt\varrho\vv+\varrho{\rm div}(\vv{\otimes}\vv)
 +(\vv\Cdot\nabla\varrho)\vv
 \\\nonumber&\qquad
 \!\!\!\stackrel{\scriptsize\eq{Euler0-fluid-accret}}{=}\!\varrho\pdt \vv-\varrho({\rm div}\vv)\vv-(\vv\Cdot\nabla\varrho)\vv
 +\varrho(\vv\Cdot\nabla)\vv+\varrho({\rm div}\vv)\vv
 +(\vv\Cdot\nabla\varrho)\vv+r_{\rm ext}^{\nothing}\vv
 \\
 &\qquad
 =\varrho\Big(\pdt \vv+(\vv\Cdot\nabla)\vv\Big)+r_{\rm ext}^{\nothing}\vv.
\label{momentum-velocity-transport-accret}\end{align}
We thus can write \eq{Euler1-fluid-accret} equivalently as
\begin{align}
\varrho\DT\vv-{\rm div}\,\DIS_{\ee}'(J,\theta;\ee(\vv))
+\nabla p=r_{\rm ext}^{\nothing}(\vv_{\rm
  ext}^{\nothing}{-}\vv)+\varrho{\bm g} -2 \varrho \oomega{\times}\vv\,.
\label{Euler1-fluid-accret-alt}
\end{align}
We now \sss multiply  of the momentum equation in the form \eq{Euler1-fluid-accret-alt} by
$\vv$, integrate over $\varOmega$, and use the Green formula, taking the
the corresponding boundary conditions into account. In particular, %\sss; by the ``test'' we mean (as usual) multiplication of the equation
%by a function in question, integrating over $\varOmega$, and using suitably Green
%formulas with relevant boundary conditions. For this particular test,
we also use 
\eee the nonhomogeneous evolution-and-transport equation
\eq{Euler0-fluid-accret} \sss multiplied \eee by $|\vv|^2/2$. We obtain 
\begin{align}
  \pdt{}\bigg(\frac\varrho2|\vv|^2\bigg)=\varrho\vv\Cdot\pdt\vv
  +\pdt\varrho\frac{|\vv|^2}2
 & =\varrho\vv\Cdot\pdt\vv
-{\rm div}(\varrho\vv)\frac{|\vv|^2}2+\frac{r_{\rm ext}^{\nothing}}2|\vv|^2\,,
\label{rate-of-kinetic-accret}\end{align}
where the last term represents the rate of increase of the kinetic energy by
the incoming mass flux $r_{\rm ext}^{}$. Thus we have
\begin{align*}
\!\!\!\int_\varOmega\varrho(\vv\Cdot\nabla)\vv\Cdot\vv\,\d\xx
  & =\int_\varOmega \varrho \vv \Cdot \nabla \left(
    \frac{|\vv|^2}{2}\right) \, \d \xx  =\int_\varOmega{\rm div}\left(\varrho \vv\frac{|\vv|^2}{2}\right)
    \, \d \xx - \int_\varOmega {\rm div} (\varrho
    \vv)\frac{|\vv|^2}{2}\, \d \xx\\
\\[-1em]&\nonumber
=-\int_\varOmega{\rm div}(\varrho\vv)\frac{|\vv|^2}2\,\d\xx
\!\!\stackrel{\scriptsize\scriptsize\eq{rate-of-kinetic-accret}}{=}\!\!\frac{\d}{\d t}
\int_\varOmega\frac\varrho2|\vv|^2\,\d\xx-\!\int_\varOmega\!
r_{\rm ext}^{\nothing}\frac{|\vv|^2}2-\varrho\pdt\vv\Cdot\vv\,\d\xx\,,
\end{align*}
i.e.,
\begin{align}\nonumber\\[-3em]
\int_\varOmega\varrho\DT\vv\Cdot\vv\,\d\xx=\frac{\d}{\d t}\int_\varOmega\frac\varrho2|\vv|^2\,\d\xx
-\int_\varOmega\frac{r_{\rm ext}^{\nothing}}2|\vv|^2\,\d\xx\,.\qquad\qquad
  \label{calculus-convective-accret}\end{align}
\ttt In order to deal with the integral of $\varrho{\bm g}\cdot \vv$,
we multiply the \eee
%Further, the gravitational force
%$\varrho{\bm g}=-\varrho\nabla\GRAVPOT$ tested by $\vv$ uses
the equation for the gravitational potential
$\GRAVPOT=4\uppi\GRAVCONST\Delta^{-1}\varrho$ in
\eq{Euler1-fluid-accret}  
by $\pdt{}\GRAVPOT$, \sss integrate on $\varOmega$, and find \eee %, which gives
\begin{align}\nonumber
\frac{\d}{\d t}\int_{\R^3}\!\!\frac{|\nabla\GRAVPOT|^2}{8\uppi\GRAVCONST}\,\d\xx
&=-\int_{\R^3}\frac{\Delta\GRAVPOT}{4\uppi\GRAVCONST}\pdt{\GRAVPOT}\,\d\xx
=-\int_\varOmega\!\!\varrho\pdt{\GRAVPOT}\,\d\xx
\\&\nonumber
=\int_\varOmega\pdt\varrho\GRAVPOT\,\d\xx
-\frac{\d}{\d t}\int_\varOmega\!\varrho\GRAVPOT\,\d\xx
\!\!\stackrel{\scriptsize\eq{Euler0-fluid-accret}}{=}\!\!\!\int_\varOmega\!\GRAVPOT\,\big(r_{\rm ext}^{\nothing}-{\rm div}(\varrho\vv)\big)\,\d\xx-\frac{\d}{\d t}\int_\varOmega\!\varrho\GRAVPOT\,\d\xx
\\&\nonumber
=\!\int_\varOmega r_{\rm ext}^{\nothing}\GRAVPOT+\varrho\vv\Cdot\nabla\GRAVPOT\,\d\xx
-\int_\varGamma\varrho\GRAVPOT\!\!\!\!\lineunder{\!\!\!\!\vv\Cdot\nn\!\!\!\!}{$=0$}\!\!\!\!\d S
-\frac{\d}{\d t}\int_\varOmega\varrho\GRAVPOT\,\d\xx\,,
\end{align}
so that
\begin{align}\nonumber\\[-2.9em]
\ttt \int_\varOmega \varrho {\bm g} \cdot \vv \, \d\xx= - \eee
  \int_\varOmega\varrho\nabla\GRAVPOT\Cdot\vv\,\d\xx= \ttt - \eee
\frac{\d}{\d t}\bigg(\int_{\R^3}\frac{|\nabla\GRAVPOT|^2}{8\uppi\GRAVCONST}\,\d\xx
+\int_\varOmega\varrho\GRAVPOT\,\d\xx\bigg)\:\ttt + \eee\!\!\int_\varOmega r_{\rm ext}^{\nothing}\GRAVPOT\,\d\xx\,.
\label{calculus-selfgravit-accret}\end{align}
Eventually, we \ttt handle \eee  the contribution from the pressure $p$ as
\begin{align}\nonumber
\nonumber
\int_\varOmega&\nabla p\Cdot\vv\,\d\xx
=\!\int_\varOmega\psi_J'(J,\theta)\,{\rm div}\,\vv\,\d\xx
=\int_\varOmega\phi'(J)\big({\rm div}\,\vv{-}v_{\rm ext}^{}\big)+v_{\rm ext}^{}\phi'(J)+\COUPLING_J'(J,\theta)\,{\rm div}\,\vv\,\d\xx
\\[-.3em]&\nonumber\!\!\!\!\stackrel{\scriptsize\eq{Euler3-fluid-selfgravit-accret}}{=}\!\!\!\int_\varOmega\phi'(J)
\frac{\DT J}J+v_{\rm ext}^{}\phi'(J)+\COUPLING_J'(J,\theta)
\,{\rm div}\,\vv\,\d\xx
\\&\nonumber=\int_\varOmega\bigg(
\Big[\frac{\phi(J)}J\Big]'+\frac{\phi(J)}{J^2}\bigg)\DT J
+v_{\rm ext}^{}\phi'(J)+\COUPLING_J'(J,\theta)\,{\rm div}\,\vv\,\d\xx
\\&\nonumber
\!\!\!\!\stackrel{\scriptsize\eq{Euler3-fluid-selfgravit-accret}}{=}\!\!\int_\varOmega
\Big[\frac{\phi(J)}J\Big]'\DT J
+\frac{\phi(J)}{J}\big({\rm div}\,\vv{-}v_{\rm ext}^{}\big)
+v_{\rm ext}^{}\phi'(J)+\COUPLING_J'(J,\theta)\,{\rm div}\,\vv\,\d\xx
\\&\nonumber
=\int_\varOmega\pdt{}\Big(\frac{\phi(J)}J\Big)
+\nabla\Big(\frac{\phi(J)}J\Big)\Cdot\vv
+\frac{\phi(J)}{J}\Big({\rm div}\,\vv{-}v_{\rm ext}^{}\Big)
+v_{\rm ext}^{}\phi'(J)+\COUPLING_J'(J,\theta)\,{\rm div}\,\vv\,
\,\d\xx
\\&\nonumber
=\frac{\d}{\d t}\!\int_\varOmega\frac{\phi(J)}J\,\d\xx
+\int_\varOmega{\rm div}\Big(\frac{\phi(J)}J\vv\Big)
+v_{\rm ext}^{}\Big(\phi'(J)-\frac{\phi(J)}{J}\Big)
+\COUPLING_J'(J,\theta)\,{\rm div}\,\vv\,\d\xx
\\&=\frac{\d}{\d t}\int_\varOmega\!\frac{\phi(J)}J\,\d\xx
+\int_\varGamma\!\frac{\phi(J)}J\hspace*{-.7em}\lineunder{\vv{\cdot}\nn}{$=0$}\hspace*{-.7em}\d S
+\int_\varOmega\!
v_{\rm ext}^{}\Big(\phi'(J)-\frac{\phi(J)}{J}\Big)
+\COUPLING_J'(J,\theta)\,{\rm div}\,\vv\,\d\xx\,,
\label{test-p-accret}\end{align}
where again $\phi(J)=\psi(J,0)$ and $\COUPLING(J,\theta):=\psi(J,\theta)-\phi(J)$.
Thus, the momentum equation \eq{Euler1-fluid-accret-alt} \sss
multiplied \eee % tested
by
$\vv$ \sss and integrated over $\varOmega$ \eee gives \begin{align}
\int_\varOmega\varrho\DT\vv\Cdot\vv+
\DIS_{\ee}'(J,\theta;\ee(\vv)){:}\ee(\vv)+\nabla p\Cdot\vv\,\d\xx
=\int_\varOmega\Big(r_{\rm ext}^{\nothing}(\vv_{\rm ext}^{\nothing}{-}\vv)-\varrho\nabla\GRAVPOT\Big)\Cdot\vv\,\d\xx\,.
\label{Euler1-fluid-accret-tested+}\end{align}
By using \eq{calculus-convective-accret}, \eq{calculus-selfgravit-accret}, and
\eq{test-p-accret}, we can rewrite \eq{Euler1-fluid-accret-tested+} into the
form of the energy-dissipation balance:
\begin{align}\nonumber
&\frac{\d}{\d t}\bigg(\int_\varOmega\hspace{-.7em}\linesunder{\frac\varrho2|\vv|^2
}{kinetic}{energy}\hspace{-.9em}
+\hspace{-.7em}\morelinesunder{\frac{\phi(J)}J+\varrho\GRAVPOT}{stored and}{gravitational}{energy}\hspace{-.7em}\d\xx
+\int_{\R^3}\hspace{-1.3em}\morelinesunder{\frac{|\nabla\GRAVPOT|^2\!\!}{8\uppi\GRAVCONST}}{energy of}{gravitational}{field}\hspace{-1.3em}\d\xx\bigg)
+\!\int_\varOmega\!\!\!\!\morelinesunder{
\,\DIS_{\ee}'(J,\theta;\ee(\vv)){:}\ee(\vv)_{_{_{_{_{_{}}}}}}\!\!}
{viscosity}{dissipation}{rate}\!\!\!\d\xx
\\[-.7em]\nonumber
&\quad=\!\int_\varOmega\!
     r_{\rm ext}^{\nothing}(\vv_{\rm ext}^{\nothing}{-}\vv)\Cdot\vv
+\frac{r_{\rm ext}^{\nothing}}2|\vv|^2+r_{\rm ext}^{\nothing}\GRAVPOT
+v_{\rm ext}^{}\Big(\frac{\phi(J)}{J}{-}\phi'(J)\Big)
-\COUPLING_J'(J,\theta)\,{\rm div}\,\vv\,\d\xx
\\&\nonumber\quad=
\!\int_\varOmega\!\bigg(\!\!\!\linesunder{\frac{r_{\rm ext}^{\nothing}}2|\vv_{\rm ext}^{\nothing}|^2
+r_{\rm ext}^{\nothing}\GRAVPOT+v_{\rm ext}^{}\Big(\frac{\phi(J)}{J}{-}\phi'(J)\Big)}
{rate of energy (power) of the material coming}{from the outer space\ \ \ \ }
\\[-1em]&\nonumber\hspace{18em}
-\!\!\!\!\!\linesunder{\COUPLING_J'(J,\theta)
\,{\rm div}\,\vv_{_{_{_{_{_{}}}}}}\!\!}{power
of}{adiabatic effects}\!\!-\!\!\!\!\!\morelinesunder{\frac{r_{\rm ext}^{\nothing}}2|\vv_{\rm ext}^{\nothing}{-}\vv|^2}
{power of friction}{of incoming}{material}
     \!\!\!\!\bigg)\,\d\xx\,.
\\[-1.7em]\label{mech-engr-selfgrav-accret}\end{align}
Summing \eq{mech-engr-selfgrav-accret} with \eq{Euler4-fluid-accret}
\sss integrated over $\varOmega$, \eee %tested by 1,
we obtain the expected total-energy balance:
\begin{align}\nonumber
\frac{\d}{\d t}&\bigg(\int_\varOmega\linesunder{\frac\varrho2|\vv|^2
}{kinetic}{energy}\hspace{-.7em}+\hspace{-.7em}\morelinesunder{\varrho\GRAVPOT}{gravita-}{tional}{energy}\hspace{-.7em}
+\hspace{-1.7em}\linesunder{\frac{\phi(J)}J+\W(J,\theta)}
{$\ \ \ \ =(\psi(J,\theta){-}\theta\psi_\theta'(J,\theta))/J$\!\!\!\!\!}{internal energy}
\hspace{-1.2em}\d\xx
+\int_{\R^3}\hspace{-1.3em}\morelinesunder{\frac{|\nabla\GRAVPOT|^2\!\!}{8\uppi\GRAVCONST}}
{energy of}{gravitational}{field}\hspace{-1em}\d\xx\bigg)
\\[-.7em]&\nonumber\qquad
=\!\int_\varOmega\!\bigg(\!\!\linesunder{\frac{r_{\rm ext}^{\nothing}}2|\vv_{\rm ext}^{\nothing}|^2+
v_{\rm ext}^{\nothing}\sss
           \Big(\frac{\psi(J,\theta) - \theta
           \psi_\theta'(J,\theta)}{J}%{-}\psi_{J}'(J,\theta)
           \Big) \eee
           }
           {rate of kinetic and \sss internal \eee %free and pressure
           energy}
{incoming from the outer space}\!\!
\\[-.5em]&\hspace*{8em}
+\!\!\morelinesunder{r_{\rm ext}^{\nothing}\GRAVPOT}{rate of energy}{of mass in}{gravitational field}+\!\!\morelinesunder{h_{\rm ext}^{\nothing}}{power of}{heat coming}{from outer space}\!\!\!
-\!\!\linesunder{\frac{r_{\rm ext}^{\nothing}}2|\vv_{\rm ext}^{\nothing}{-}\vv|^2}
{power of friction}{of incoming material}\!\!\!\!\bigg)\d\xx\,.
\label{total-engr-selfgrav-accret}\end{align}
It should be noted that  we are assuming that the temperature of
the incoming material corresponds to that of the particular point in
space time, where the material occurs. On the other hand, the extra heat is
accounted for by $h_{\rm ext}^{\nothing}$.

Notably, as we started from   \eq{entropy-eq+}, the model complies with the
2nd law of thermodynamics. More specifically, from \eq{entropy-eq+}, we have
the entropy balance
\begin{align}\label{ent-balance}
& \hspace*{-2em}\frac{\d}{\d t}\!\!\!\!\linesunder{\int_\varOmega\eta\,\d\xx}{total}{entropy}
\!\!\!\!=\int_\varOmega\!\!\!\!\lineunder{\frac{\kappa|\nabla\theta|^2\!}{\theta^2}
  +\frac{\xi+h_{\rm ext}^{}}{\theta}+v_{\rm ext}^{}\eta}{entropy production rate}\!\!\!\d\xx\,.
\end{align}
In particular, when $\xi\ge0$, $h_{\rm ext}^{}\ge0$, and $v_{\rm ext}^{}\ge0$,
the total entropy is nondecreasing in time.

\begin{remark}[{\sl Friction of incoming material}]\upshape
The last nonpositive terms in \eq{mech-engr-selfgrav-accret} and
\eq{total-engr-selfgrav-accret} arise from friction when the incoming material
has a different velocity $\vv_{\rm ext}^{}$ than the current velocity of the
medium $\vv$. They would vanish if we would have imposed 
that these velocities coincide, consistently with the assumption on 
temperature and on volume compression (determined by $J$).
This latter point is however quite acceptable because the incoming material
is presumably very cold (i.e.\ $\theta\sim0$) and very diluted ($J>\!\!>1$),
similarly as the material in  those   regions of $\varOmega$ far from the accreting
planet. Yet, there is no reason why
$\vv_{\rm ext}^{\nothing}-\vv$ should be small, which is why we kept the
general discrepancy between them in the model.
\end{remark}

\begin{remark}[{\sl Outflowing material}]\label{rem-outflowing}\upshape
In reality, not all material which enters the vicinity of the accreting planet
falls on the planet. Some portion of the  incoming  material leaves again
to the outer space, as indicated in Figure~\ref{fig-accretion}. This suggests
that $r_{\rm ext}^{}$ may also be negative. Of course, it should be dependent on
the mass density and the current velocity, not to withdraw material if there
is none (which would nonphysically lead to  a negative density  $\varrho$) or
if its velocity vector is not pointing towards the outer space. Thus, if
$r_{\rm ext}^{}(t,\xx)<0$, one should rather consider
$r_{\rm ext}^{}=r_{\rm ext}^{}(\xx,\varrho,\vv)$. Also, it is reasonable to consider
the  outflowing   material having the velocity $\vv$, so
$\vv_{\rm ext}^{}(t,\xx)=\vv(t,\xx)$.
\end{remark}

\section{A two-component metal-silicate model}\label{sect-two-comp}
%        ~~~~~~~~~~~~~~~~~~~~~~~~~~~~~~~~~~~~

\def\S{_\text{\sc s}^{}}
\def\M{_\text{\sc m}^{}}
\def\Sxi{_\text{\sc s}^{\bm\xi\S}}
\def\Mxi{_\text{\sc m}^{\bm\xi\M}}
\def\Sref{_\text{\sc s,{\rm ref}}^{}}
\def\Mref{_\text{\sc m,{\rm ref}}^{}}
\def\Srefxi{_\text{\sc s,{\rm ref}}^{\bm\xi\S}}
\def\Mrefxi{_\text{\sc m,{\rm ref}}^{\bm\xi\M}}
\def\Sext{_\text{\sc s,{\rm ext}}^{}}
\def\Mext{_\text{\sc m,{\rm ext}}^{}}
\def\Sextxi{_\text{\sc s,{\rm ext}}^{}}
\def\Mextxi{_\text{\sc m,{\rm ext}}^{}}

We now extend the above single-material model  in order to model
differentiation driven by self-gravitation. To this aim, we necessarily need to
consider a multi-component system, with single constituents having
different densities and different thermo-viscoelastic
responses. A minimal scenario is  that of  considering two
components,  namely, {\it metal}
and {\it silicate}.  This is actually  the coarsest description of a usual
structure of terrestrial type planets,  which have  a heavier metal core and 
a lighter silicate mantle.  We hence  will distinguish these components by the
indexes ``M'' and ``S''.

Both components have their own velocity $\vv\S$ and $\vv\M$ and
temperature $\theta\S$ and $\theta\M$, respectively.
In addition, each component has its own density $\varrho\S$ and
$\varrho\M$ and a viscoelastic response governed by  the 
referential free energies
$\psi\S=\psi\S(J,\theta)$ and $\psi\M=\psi\M(J,\theta)$ and  the 
dissipation potentials $\DIS\M=\DIS\M(J,\theta;\ee)$ and $\DIS\S=\DIS\S(J,\theta;\ee)$.
We again consider the split
\begin{subequations}\begin{align}
&\psi\S(J,\theta)=\phi\S(J)+\COUPLING\S(J,\theta)\ \ \ \ \text{ with }\ \ \ \ 
\COUPLING\S(J,0)=0\,,\ \  \ \text{ and}
\\
&\psi\M(J,\theta)=\phi\M(J)+\COUPLING\M(J,\theta)\ \ \ \text{ with \ \ }\ 
\COUPLING\M(J,0)=0\,.
\end{align}\end{subequations}

Both components enter our open system with their own mass density rates
$r\Mext=r\Mext(t,\xx)$ and $r\Sext=r\Sext(t,\xx)$,
volumes $v\Mext=v\Mext(t,\xx)$ and $v\Sext=v\Sext(t,\xx)$, 
velocities $\vv\Mext=\vv\Mext(t,\xx)$ and $\vv\Sext=\vv\Sext(t,\xx)$, and
heat rates  $h\Mext=h\Mext(t,\xx)$ and $h\Sext=h\Sext(t,\xx)$.
 We impose the relation \eqref{r-vs-v} on both components, namely,
$$
r\Mext=\frac{\varrho_{0,{\sc m}}v\Mext}{J\Mext}\quad\text{and}\quad
r\Sext=\frac{\varrho_{0,{\sc s}}v\Sext}{J\Sext}\,.
$$

It is important to model the interaction between these two components. First,
when moving with different velocities, one should expect some friction-like
force between them. Let us consider it (for simplicity) as linearly dependent
on the difference of velocities through a friction coefficient
$f=f(\varrho\M,\varrho\S)$. Second, when $\theta\S\ne\theta\M$, one should
consider an exchange of heat between these two components, depending through
a heat-exchange coefficient $k=k(\varrho\M,\varrho\S)$ (for simplicity
linearly) on the difference of the temperatures. Third, an important modelling
ingredient is a phase separation in case of high densities.
We will promote such separation by introducing a suitable mixing energy
$\phi_{\rm mix}=\phi_{\rm mix}(J\M,J\S)$\sss, modeling the tendency of
the components (here, metals and silicates) to separate unless mutually
diluted. \eee
%reflecting the philosophy that
%the particular components (i.e.\ here metals and silicates) do not want
%to co-exist on the same spots except if they are only diluted\eee.
In contrast to $\psi\M$ and $\psi\S$,  the  mixing energy $\phi_{\rm mix}$ is
related to different volumes of metals and silicates  and is hence to be
considered with respect to actual volumes. To be more specific, let us
assume  that such mixing energy has the form 
\begin{align}
\phi_{\rm mix}(J\M,J\S)=\varkappa\begin{cases}
(J\M J\S)^{-\alpha}+\alpha J\M J\S-\alpha-1&\text{if }\ J\M J\S<1\,,
\\\qquad\qquad0&\text{if }\ J\M J\S\ge1
\end{cases}
\label{example}\end{align}
with some  given and large $\varkappa>0$ (in physical unit Pa\,=\,J/m$^3$)
and some $\alpha>0$. Note that $\phi_{\rm mix}$ in \eq{example} is continuously
differentiable. \sss We refer to
\cite[Sects 5.2 and 5.4]{PoPrOC01PGL} for a discussion of different
possible mixing-energy forms, as well as for references. \eee %various mixing and combining rules
%and many references there.\eee 

The self-gravitating system \eq{Euler-fluid-accret} is then considered doubled
to describe the metallic and the silicate components separately, and augmented
by the mentioned interaction mechanisms. Altogether, it yields the system
\begin{subequations}\label{Euler-fluids-accret}\begin{align}
&
\pdt{\varrho\M}+{\rm div}(\varrho\M\vv\M)=r\Mextxi\ \ \ \text{ and }\ \ \ 
\label{Euler0-fluids-accret}
\pdt{\varrho\S}+\,{\rm div}(\varrho\S\vv\S)=r\Sextxi\,,\\
\nonumber
&\pdt{}(\varrho\M\vv\M)+{\rm div}\Big(\varrho\M\vv\M{\otimes}\vv\M)-
[\DIS\M]_{\ee}'(J\M,\theta\M;\ee(\vv\M))
\Big)
+\nabla p\M
\\[-.4em]&\nonumber\hspace{4em}
=r\Mextxi\vv\Mextxi
+f(\varrho\M,\varrho\S)(\vv\M{-}\vv\S)
+\big[\phi_{\rm mix}\big]_{J\M}'\!(J\M,J\S)\nabla J\M+\varrho\M{\bm g}
-2\varrho\M\oomega{\times}\vv\M\ \text{ and}
\\\nonumber
&\pdt{}(\varrho\S\vv\S)+{\rm div}\Big(\varrho\S\vv\S{\otimes}\vv\S)-
[\DIS\S]_{\ee}'(J\S,\theta\S;\ee(\vv\S))\Big)
+\nabla p\S
\\[-.4em]&\nonumber\hspace{4em}
=r\Sextxi\vv\Sextxi+f(\varrho\M,\varrho\S)(\vv\S{-}\vv\M)
 +\big[\phi_{\rm mix}\big]'_{J\S}\!(J\M,J\S)\nabla J\S+\varrho\S{\bm
           g} -2\varrho\S\oomega{\times}\vv\S\ \text{ with}
\\&\nonumber\hspace{4em}
p\M=-\def\Sextxi{_\text{\sc s,{\rm ext}}^{\bm\xi\S}}
\big[ \psi\M\big]'_{ J\M}(J\M,\theta\M)
-{ J\M}\big[\phi_{\rm mix}\big]'_{J\M}\!(J\M,J\S)\,,
\\&\nonumber\hspace{4em}
p\S=-\big[\psi\S\big]'_{ J\S}(J\S,\theta\S)
-{ J\S}\big[\phi_{\rm mix}\big]_{J\S}'\!(J\M,J\S)\,,
\\&\hspace{4em}
{\bm g}(\cdot)=\GRAVCONST\!\int_\varOmega\frac{(\varrho\M(\widetilde\xx){+}\varrho\S(\widetilde\xx))(\widetilde\xx-\cdot)\!}{|\widetilde\xx-\cdot|^3}\,\d\widetilde\xx\,,
\label{Euler1-fluids-accret}
\\\nonumber
&\pdt{\W\M}+{\rm div}\Big(\vv\M\W\M{-}\kappa\M(J\M,\theta\M)\nabla\theta\M\Big)
=\big[\DIS\M\big]_{\ee}'(J\S,\theta\M;\ee(\vv\M)){:}\ee(\vv\M)
\\[-.4em]&\nonumber\hspace{4em}
+\big(\big[\psi\M\big]_J'(J\M,\theta\M){-}\phi_\text{\sc m}'(J\M)\big)
\big({\rm div}\,\vv\M{-}v\Mextxi\big)
\sss+\eee
% \frac{\COUPLING\M(J\M,\theta\M)}{J\M}v\Mextxi\COMMENT{\ OR\
           \sss \W\M v\Mextxi \eee %?????}
\\[-.4em]&\nonumber\hspace{4em}
+k(\varrho\M,\varrho\S)(\theta\S{-}\theta\M)
+\frac12f(\varrho\M,\varrho\S)|\vv\M{-}\vv\S|^2\!+h\Mextxi\,,
\\[-.2em]&\nonumber\hspace{16em}\text{with }\ 
\W\M=\frac{\COUPLING\M(J\M,\theta\M){-}\theta\M[\COUPLING\M]_\theta'(J\M,\theta\M)}{J\M}
\ \ \ \text{ and }
\\\nonumber
&\pdt{\W\S}+{\rm div}\Big(\vv\S\W\S{-}
\kappa\S(J\S,\theta\S)\nabla\theta\S\Big)
=\big[\DIS\S\big]_{\ee}'(J\S,\theta\S;\ee(\vv\S)){:}\ee(\vv\S)
\\[-.4em]&\nonumber\hspace{4em}
+\big(\big[\psi\S\big]_J'(J\S,\theta\S){-}\phi_\text{\sc s}'(J\S)\big)
\big({\rm div}\,\vv\S{-}v\Sextxi\big)
           \sss + %\frac{\COUPLING\S(J\S,\theta\S)}{J\S}v\Sextxi\COMMENT{\ OR\
           \W\S v\Sextxi \eee
\\[-.4em]&\nonumber\hspace{4em}
+k(\varrho\M,\varrho\S)(\theta\M{-}\theta\S)
+\frac12f(\varrho\M,\varrho\S)|\vv\M{-}\vv\S|^2+h\Sextxi
\\[-.2em]&\hspace{16em}\text{with }\ 
\
\label{Euler4-fluids-accret}
\W\S\!=\frac{\COUPLING\S(J\S,\theta\S){-}\theta\S[\COUPLING\S]_\theta'(J\S,\theta\S)
}{J\S}\,,\!\!
\\[-.2em]
&
\pdt{J\M}=({\rm div}\,\vv\M-v\Mextxi)J\M-\vv\M\Cdot\nabla J\M\ \ \text{ and }\ \ 
\pdt{J\S}=({\rm div}\,\vv\S-v\Sextxi)J\S-\vv\S\Cdot\nabla J\S\,.
\label{Euler3-fluids-selfgravit-accret}
\end{align}\end{subequations}

As in Section~\ref{sec-single}, to establish the balance of mass, momentum, 
energy, and entropy, we must prescribe boundary
conditions. We again consider the  boundary $\varGamma$ of the 
domain $\varOmega$  to be impermeable and friction-free and
$\varOmega$ to be thermally insulated from outer space by prescribing
\begin{align}\nonumber
\vv\M\Cdot\nn=\vv\S\Cdot\nn=0\,,\ \ \ \ 
&[\DIS_{\ee}'(J\M,\theta\M;\ee(\vv\M))\nn]_\text{\sc t}=0\,,\ \ \ 
\\&[\DIS_{\ee}'(J\S,\theta\S;\ee(\vv\S))\nn]_\text{\sc t}=0\,,
\ \text{ and}\ \nabla\theta\M\Cdot\nn=\nabla\theta\S\Cdot\nn=0
\ \ \text{ on }\ \varGamma\,.
\label{Euler-fluids-accret-BC}\end{align}

The mass balance now contains source terms due to occurrence of incoming
silicates and metals by the rates $r\Mref$ and $r\Sref$, respectively. The
{\it total-mass balance} can be obtained by integrating both equations in
\eq{Euler0-fluids-accret} over $\varOmega$:
\begin{align}
\frac{\d}{\d t}\!\!\!\!\linesunder{\int_\varOmega\varrho\M{+}\varrho\S\,\d\xx}
{total mass in}{the domain $\varOmega$}
\!\!\!\!=\!\!\!\!\linesunder{\int_\varOmega r\Mextxi+r\Sextxi\,\d\xx}{incoming}{mass rate}\!\!.
\end{align}

The momentum balance now contains source terms due to  the
additional  momenta
of both silicates and metals with the rates  $r\Mextxi\vv\Mextxi$ and
$r\Sextxi\vv\Sextxi$, respectively. It can be obtained by
integrating both equations in \eq{Euler1-fluids-accret} over $\varOmega$.
The mixing forces merge by Stokes' theorem as
\begin{align}\nonumber
\int_\varOmega\big[\phi_{\rm mix}\big]_{J\M}'\!\!(J\M,J\S)\nabla J\M
&+\big[\phi_{\rm mix}\big]_{J\S}'\!(J\M,J\S)\nabla J\S\,\d\xx
\\[-.4em]&=\int_\varOmega\nabla\big[\phi_{\rm mix}(J\M,J\S)\big]\,\d\xx
=\int_\varGamma\phi_{\rm mix}(J\M,J\S)\nn\,\d S\,.
\label{mix-calculus}\end{align}
Consistently with  the invariance  of \eq{Euler-fluids-accret} with respect to
additive constants to $\phi_{\rm mix}$, also \eq{mix-calculus} is invariant
since $\int_\varGamma\nn\,\d S=\bm0$. The frictional forces cancel each other
pointwise. As for the self-gravitational momenta, they cancel in overall sum
by the argument \eq{self-graviating-momentum} now  applied to 
$\varrho=\varrho\M{+}\varrho\S$.  On the other hand, the Coriolis contributions
are added.  Altogether, the {\it total-momentum balance} reads as
\begin{align}
\frac{\d}{\d t}\!\!\!\!\lineunder{\int_\varOmega\varrho\M\vv\M+\varrho\S\vv\S\,\d\xx}{total momentum}
\!\!\!\!&=\!\!\!\!\lineunder{\int_\varOmega
r\Mextxi\vv\Mextxi+r\Sextxi\vv\Sextxi\,\d\xx}{incoming momentum rate}\!\!
\nonumber
  \\&
  +\!\!\!\!\!\!\lineunder{\int_\varGamma\phi_{\rm mix}(J\M,J\S)\nn\,\d
  S}{incoming mixing energy rate}\!\!\!\!\!\!\!-\!\!\!\!
  \lineunder{ \int_\varOmega 2\varrho\M\oomega{\times}\vv\M\,\d\xx
  +2\varrho\S\oomega{\times}\vv\S\,\d\xx}{ the Coriolis
  momentum rate}\!\!\!.
\label{mometum-tot}\end{align}

The next task is to specify the {\it energy-dissipation balance}. To this aim,
the individual momentum equations in  \eq{Euler1-fluids-accret} is to be
\sss multiplied \eee by $\vv\M$ and $\vv\S$,   respectively, \sss and then integrated on
$\varOmega$. \eee

For the  metal component,  as in \eq{momentum-velocity-transport-accret}, we
have
\begin{align}
&\pdt{}(\varrho\M\vv\M)+{\rm div}(\varrho\M\vv\M{\otimes}\vv\M)
=\varrho\M\Big(\pdt{\vv\M}+(\vv\M\Cdot\nabla)\vv\M\Big)+r\Mextxi\vv\M\,.
\label{momentum-velocity-transport-accret+}\end{align}
The analogous computation holds also for the silicate component. For the power
of the gravitational forces, i.e., of $\varrho\M{\bm g}$
\sss multiplied \eee by $\vv\M$ and of $\varrho\S{\bm g}$ \sss multiplied \eee by $\vv\S$ with
${\bm g}=-\GRAVCONST\nabla\GRAVPOT$, as in \eq{calculus-selfgravit-accret},
we have
\begin{align}
\!\!\int_\varOmega\!
\big(\varrho\M\vv\M{+}\varrho\S\vv\S\big)\Cdot{\bm g}\,\d\xx=\!
\int_\varOmega\!\big(r\Mextxi{+}r\Sextxi\big)\GRAVPOT\,\d\xx-
\frac{\d}{\d t}\bigg(\int_{\R^3}\!\!\frac{|\nabla\GRAVPOT|^2\!\!}{8\uppi\GRAVCONST}\,\d\xx
+\!\!\int_\varOmega\!(\varrho\M{+}\varrho\S)\GRAVPOT\,\d\xx\bigg).\!
\label{calculus-selfgravit-accret+}\end{align}

\sss Forces \eee and the pressure gradient arising
from the mixing energy $\phi_{\rm mix}$ in the momenta equations
\eq{Euler1-fluids-accret} give
\begin{align}\nonumber
&-\!\int_\varOmega\!\bigg(\big[\phi_{\rm mix}\big]_{J\M}'(J\M,J\S)\nabla J\M
+\nabla\Big(J\M\big[\phi_{\rm mix}\big]_{J\M}'\!(J\M,J\S)\Big)\bigg)\Cdot\vv\M\,\d\xx
\\[-.4em]&\nonumber\qquad\qquad\qquad\qquad
-\!\int_\varOmega\!\bigg(\big[\phi_{\rm mix}\big]_{J\S}'(J\M,J\S)\nabla J\S
+\nabla\Big(J\S\big[\phi_{\rm mix}\big]_{J\S}'\!(J\M,J\S)\Big)\bigg)\Cdot\vv\S\,\d\xx
\\[-.4em]&\nonumber\quad
\!\!\!\!\stackrel{\scriptsize\eq{Euler3-fluids-selfgravit-accret}}{=}\!\!\!
\int_\varOmega\!\bigg(\big[\phi_{\rm mix}\big]_{J\M}'\!(J\M,J\S)\,\big(J\M{\rm div}\,\vv\M-
\vv\M\Cdot\nabla J\M\big)
\\[-.7em]&\nonumber\qquad\qquad\qquad\qquad
+\big[\phi_{\rm mix}\big]_{J\S}'\!(J\M,J\S)\,\big(J\S{\rm div}\,\vv\S
-\vv\S\Cdot\nabla J\S\big)\bigg)\,\d\xx
\\&\nonumber\quad
=\int_\varOmega\!\bigg(\big[\phi_{\rm mix}\big]_{J\M}'(J\M,J\S)\,\pdt{J\M}
+\big[\phi_{\rm mix}\big]_{J\S}'(J\M,J\S)\,\pdt{J\S}
\\[-.7em]&\nonumber\qquad\qquad\qquad\qquad
+J\M\big[\phi_{\rm mix}\big]_{J\M}'\!(J\M,J\S)v\Mextxi
+J\S\big[\phi_{\rm mix}\big]_{J\S}'\!(J\M,J\S)v\Sextxi\bigg)\,\d\xx
\\&\quad
=\frac{\d}{\d t}\int_\varOmega\!\phi_{\rm mix}(J\M,J\S)\,\d\xx
+\!\int_\varOmega\!J\M\big[\phi_{\rm mix}\big]_{J\M}'\!(J\M,J\S)v\Mextxi\!
+J\S\big[\phi_{\rm mix}\big]_{J\S}'\!(J\M,J\S)v\Sextxi\,\d\xx\,.
\label{mixing-calculus}
\end{align}

Eventually, we should evaluate also the contribution from the incoming free
energy. Then the calculation \eq{test-p-accret} exploits
\eq{Euler3-fluids-selfgravit-accret}
and modifies as
\begin{align}\nonumber
\nonumber
-\!\int_\varOmega&\nabla\big[\psi\M\big]_{ J\M}'(J\M,\theta\M)\Cdot\vv\M\,\d\xx
=\frac{\d}{\d t}\int_\varOmega\!\frac{\phi\M(J\M)}{J\M}\,\d\xx
\\&
+\!\int_\varOmega\!
v\Mextxi\Big({\phi'\!}\M(J\M)-\frac{\phi\M(J\M)}{J\M}\Big)
+  [\COUPLING\M]_{J\M}' (J\M,\theta\M)\,{\rm div}\,\vv\M\,\d\xx
\label{test-pM-accret}\end{align}
with $\COUPLING\M(J,\theta)=\psi\M(J,\theta)-\phi\M(J)$. Analogous calculations
hold also for the silicate pressure contribution $\psi_J'(J\S,\theta\S)$. The
above example \eq{example} gives the pressures
$-J\M[\phi_{\rm mix}]_{J\M}'\!(J\M,J\S)
= -\alpha\varkappa  J\M((J\M)^{-a-1}(J\S)^{-a}-J\S)
=\alpha\varkappa   J\M J\S(1{-}1/(J\M J\S)^{a+1})$ and
 equally for  $-J\S\big[\phi_{\rm mix}\big]_{J\S}'\!(J\M,J\S)$
whenever $J\M J\S< 1$,  otherwise for $J\M J\S\ge1$ these pressures vanish.

Summing up and merging with \eq{mixing-calculus}, we obtain
\begin{align}\nonumber
\int_\varOmega&\nabla p\M\Cdot\vv\M+\nabla p\S\Cdot\vv\S
+\big[\phi_{\rm mix}\big]_{J\M}'(J\M,J\S)\nabla J\M
+\big[\phi_{\rm mix}\big]_{J\S}'(J\M,J\S)\nabla J\S
\,\d\xx
\\&\nonumber
=\frac{\d}{\d t}\int_\varOmega\!\frac{\phi\M(J\M)}{J\M}+\frac{\phi\M(J\S)}{J\S}
+\phi_{\rm mix}(J\M,J\S)\,\d\xx
\\&\nonumber
+\int_\varOmega\bigg(v\Mextxi\Big({\phi'\!}\M(J\M)-\frac{\phi\M(J\M)}{J\M}
+J\M\big[\phi_{\rm mix}\big]_{J\M}'\!(J\M,J\S)\Big)
+  [\COUPLING\M]_{J\M}' (J\M,\theta\M)\,{\rm div}\,\vv\M
\\[-.2em]&\hspace{3.3em}
+v\Sextxi\Big({\phi'\!}\S(J\S)-\frac{\phi\S(J\S)}{J\S}
+J\S\big[\phi_{\rm mix}\big]_{J\S}'\!(J\M,J\S)\Big)
+  [\COUPLING\S]_{J\S}' (J\S,\theta\S)\,{\rm div}\,\vv\S\bigg)\,\d\xx\,.
\label{pressure-energy-mix}\end{align}

Furthermore, we note that the friction forces in both equations, when
\sss multiplied by the respective velocities and integrated over
$\varOmega$  \eee
as above, merge and yield
\begin{align}
f(\varrho\M,\varrho\S)(\vv\M{-}\vv\S)\Cdot\vv\M
+f(\varrho\M,\varrho\S)(\vv\S{-}\vv\M)\Cdot\vv\S
=f(\varrho\M,\varrho\S)|\vv\M{-}\vv\S|^2\,,
\label{friction-energy-mix}\end{align}
which gives the frictional dissipation rate. This contributes as a heat source
equally to both heat equations in \eq{Euler4-fluids-accret}.

By using  \eq{momentum-velocity-transport-accret+}-\eq{mixing-calculus},
\eq{pressure-energy-mix}-\eq{friction-energy-mix},
we  obtain the energy-dissipation balance: 
\begin{align}\nonumber
&\!\!\frac{\d}{\d t}\bigg(\int_\varOmega\hspace{-.8em}
\linesunder{\frac{\varrho\M}2|\vv\M|^2
+\frac{\phi\M(J\M)}{J\M}}{kinetic and stored}{energy of metals}\hspace{-.8em}
+\hspace{-.8em}\linesunder{\frac{\varrho\S}2|\vv\S|^2
+\frac{\phi\S(J\S)}{J\S}}{kinetic and stored}{energy of silicates}\hspace{-.8em}
+\hspace{-.8em}\linesunder{(\varrho\M{+}\varrho\S)\GRAVPOT}
{gravitational}{energy}\hspace{-.8em}
+\hspace{-.8em}\linesunder{\phi_{\rm mix}(J\M,J\S)}{mixing}{energy}\hspace{-.8em}
\,\d\xx
+\!\int_{\R^3}\hspace{-1.5em}\morelinesunder{\frac{|\nabla\GRAVPOT|^2\!\!}{8\uppi\GRAVCONST}}{energy of}{gravitational}{field}\hspace{-1.3em}\d\xx\!\bigg)
\\[-.7em]&\nonumber
+\!\int_\varOmega\hspace{-.8em}\morelinesunder{\DIS_{\ee}'(J\M,\theta\M;\ee(\vv\M)){:}\ee(\vv\M)_{_{_{_{_{_{}}}}}}\!\!}
{metal viscosity}{dissipation}{rate}\hspace{-.8em}+\hspace{-.8em}\morelinesunder{\DIS_{\ee}'(J\S,\theta\S;\ee(\vv\S)){:}\ee(\vv\S)_{_{_{_{_{_{}}}}}}\!\!}
{silicate viscosity}{dissipation}{rate}\hspace{-.8em}
+\hspace{-.8em}\morelinesunder{f(\varrho\M,\varrho\S)|\vv\M{-}\vv\S)|^2}{dissipation rate due to}{friction between}{metals and silicates}\hspace{-.8em}\d\xx
\\&\nonumber\quad=
\!\int_\varOmega\!\bigg(
\hspace{-.8em}\linesunder{r\Mextxi\Big(\frac{|\vv\Mextxi|^2\!\!}2
{+}\GRAVPOT\Big)+v\Mextxi\Big(\frac{\phi\M(J\M)}{J\M}{-}\phi{\M}'(J\M)
{ -J\M\big[\phi_{\rm mix}\big]_{J\M}'\!(J\M,J\S)}\Big)\!\!
}{energy rate of metals coming}{from the outer space}\!\!\!
\\&\nonumber\hspace{4.5em}
+\hspace{-.8em}\linesunder{r\Sextxi\Big(\frac{|\vv\Sextxi|^2\!\!}2+\GRAVPOT\Big)
+v\Sextxi\Big(\frac{\phi\S(J\S)}{J\S}{-}\phi{\S}'(J\S)
{ -J\S\big[\phi_{\rm mix}\big]_{J\S}'\!(J\M,J\S)}\Big)\!\!
}
{energy rate of silicates coming}{from the outer space}\!\!\!\!
\\[-.3em]&\hspace{1em}
-\!\!\!\!\!\linesunder{[\COUPLING\M]_{J\M}' (J\M,\theta\M)
\,{\rm div}\,{\vv\M}_{_{_{_{_{_{}}}}}}\!\!}{power of adiabatic}{effects in metals}\!\!\!\!
-\!\!\!\!\!\linesunder{[\COUPLING\S]_{J\S}' (J\S,\theta\S)
\,{\rm div}\,{\vv\S}_{_{_{_{_{_{}}}}}}\!\!}{power of adiabatic}{effects in silicates}\!\!\!\!
-\hspace{-.8em}\linesunder{\frac{r\Mextxi}2|\vv\Mextxi{-}\vv\M|^2}
{power of friction}{of incoming metals}\hspace{-.8em}
-\hspace{-.8em}\linesunder{\frac{r\Sextxi}2|\vv\Sextxi{-}\vv\S|^2}
{power of friction}{of incoming silicates}\hspace{-.8em}\bigg)\d\xx\,.
\label{mech-engr-selfgrav-accret-mix}\end{align}
Summing \eq{mech-engr-selfgrav-accret-mix} with both equations in
\eq{Euler4-fluids-accret} \sss integrated over $\varOmega$ \eee and realizing that the mutual exchange
of heat by the difference $\theta\M-\theta\S$ cancels,
we obtain the expected total-energy balance:%\COMMENT{\ TO\ CHECK:}
\begin{align}\nonumber
\frac{\d}{\d t}&\int_\varOmega\bigg(\hspace{-.8em}\linesunder{\frac{\varrho\M}2|\vv\M|^2\!
+\frac{\phi\M(J\M)}{J\M}+\W\M(J\M,\theta\M)+\varrho\M\GRAVPOT}
{kinetic, internal, and gravitational}{energy of metals}\hspace{-.8em}
+\hspace{-.8em}\linesunder{\frac{\varrho\S}2|\vv\S|^2\!
+\frac{\phi\S(J\S)}{J\S}+\W\S(J\S,\theta\S)+\varrho\S\GRAVPOT}
{kinetic, internal, and gravitational}{energy of silicates}\hspace{-.8em}
\\[-.0em]&\nonumber
+\hspace{-.8em}\linesunder{\phi_{\rm mix}(J\M,J\S)}{mixing}{energy}\hspace{-.8em}
\bigg)\d\xx+\frac{\d}{\d t}\int_{\R^3}\hspace{-1.2em}\morelinesunder{\frac{|\nabla\GRAVPOT|^2\!\!}{8\uppi\GRAVCONST}}
{energy of}{gravitational}{field}\hspace{-1.2em}\d\xx\bigg)
=\!\int_\varOmega\!\bigg(\hspace{-.8em}\morelinesunder{(r\Mextxi{+}r\Sextxi)\GRAVPOT}{rate of energy}{of mass in}{gravitational field}\hspace{-.8em}+\hspace{-.8em}\morelinesunder{h\Mextxi{+}h\Sextxi}{power of}{heat coming}{from outer space}\hspace{-.8em}
\\[-.3em]&\nonumber\qquad
+\hspace{-.8em}\linesunder{\frac{r\Mextxi}2|\vv\Mextxi|^2
+v\Mextxi\Big(\sss \frac{\psi\M(J\M,\theta\M)-\theta\M[\psi\M]_{J\M}' (J\M,\theta\M)}{J\M} \eee {+}{-J\M\big[\phi_{\rm mix}\big]_{J\M}'\!(J\M,J\S)}\Big)}
           {rate of kinetic and internal energy of}{metals coming from the outer space}\hspace{-.8em}
\\[-.2em]&\nonumber\qquad
+\hspace{-.8em}\linesunder{\frac{r\Sextxi}2|\vv\Sextxi|^2
+v\Sextxi\Big(\sss \frac{\psi\S(J\S,\theta\S)-\theta\S[\psi\S]_{J\S}' (J\S,\theta\S)}{J\S} \eee-J\S\big[\phi_{\rm mix}\big]_{J\S}'\!(J\M,J\S)\Big)}
           {rate of kinetic and internal energy of}{silicates coming from the outer space}\hspace{-.8em}
\\[-.2em]&\hspace{14em}
-\hspace{-.8em}\linesunder{\frac{r\Mextxi}2|\vv\Mextxi{-}\vv\M|^2}
{power of friction}{of incoming metals}\hspace{-.8em}
-\hspace{-.8em}\linesunder{\frac{r\Sextxi}2|\vv\Sextxi{-}\vv\S|^2}
{power of friction}{of incoming silicates}\hspace{-.8em}\bigg)\d\xx\,.
\label{total-engr-selfgrav-accret-mix}\end{align}

\begin{remark}[{\sl Momentum balance \eq{mometum-tot}}]\upshape
The last boundary-integral terms in \eq{mix-calculus} and 
\eq{mometum-tot} do not vanish in general and may possibly corrupt the expected
total-momentum balance. If the mixing forces were not involved in
\eq{Euler1-fluids-accret}, the total-momentum balance would hold. On the other
hand,  in this case  the energetics
would be corrupted, cf.\ the calculus \eq{mixing-calculus}. This inconsistency
is the price  to pay   for the simplified mixture theory adopted in
this section. Yet, in the 
applicable setting  these boundary terms in \eq{mix-calculus}
and \eq{mometum-tot} are small or even zero.
 One can assume  that the material will be arbitrarily diluted far
from the accreting planet, in particular on $\varGamma$, inducing
that both $J\M$ and $J\S$ are very large on $\varGamma$. In the
case of Example \eq{example}, this makes $\phi_{\rm mix}^{}$ vanishing
and the desired momentum balance is satisfied.
\end{remark}

\begin{remark}[{\sl Entropy imbalance}]\upshape
The exchange of  heat  between  the  two components
in \eq{Euler4-fluids-accret} brings 
an additional entropy production rate into the system.  More
precisely, 
\eq{ent-balance} written in a multi-component way for
$\eta=\eta\M+\eta\S$  features,  beside
the frictional entropy source $f(\varrho\M,\varrho\S)
|\vv\M{-}\vv\S|^2(\theta\M{+}\theta\S)/(\theta\M\theta\S)\ge0$,
also by the entropy production 
$$
k(\varrho\M,\varrho\S)\frac{\theta\S{-}\theta\M}{\theta\M}
+k(\varrho\M,\varrho\S)\frac{\theta\M{-}\theta\S}{\theta\S}
=k(\varrho\M,\varrho\S)\frac{|\theta\M{-}\theta\S|^2}{\theta\M\theta\S}\ge0\,.
$$
\end{remark}

\section{Formal stability and analytical remarks}\label{sec-anal}
%        ~~~~~~~~~~~~~~~~~~~~~~~~~~~~~~~~~~~~~~~~~~~~~~~~~~~

Providing a rigorous  analytical justification of the above models is
highly nontrivial and is beyond the scope of this paper. In this section,
we however present some formal stability argument, which could be the
departing point  for building computational amenable numerical strategies.

The analysis of the two-component model in Section~\ref{sect-two-comp}
 does not essentially differ from the one of the single-component
model. In fact, the only added term in the  energetics
\eq{mech-engr-selfgrav-accret-mix}--\eq{total-engr-selfgrav-accret-mix}
 of the two-component system  is $\phi_{\rm mix}(J\M,J\S)$
 which does not corrupt the stability argument if  $\phi_{\rm mix}\ge0$.
 We hence resort in analyzing the single-component model only.  
Beside \eq{r-vs-v}, several  assumptions are needed for the data $\psi$, $\DIS$,
$v_{\rm ext}^{}$ (now considered depending also on $J$), and $h_{\rm ext}^{}$, as
well as on the initial conditions.  In the following, we use the
symbols 
$C(\cdot)$ for the space of continuous functions and $L^q(\cdot)$
the space of function whose $q$-power is Lebesgue integrable. 
 We shall assume the following 
\begin{subequations}\label{Euler-fluid-accret-ass}
\begin{align}
 & \exists\epsilon>0,\ K>0, \ \beta>1, \ q \in (1,+\infty) \
    \ \text{such that}\nonumber\\
& \label{compressibility}
\psi\in C^1((0,+\infty){\times}\R^+),\ \ 
\liminf_{J\to0+}J^{1+\epsilon}\phi(J)>0\ \ \text{ with }\ \phi(J)=\psi(J,0)\,,
\\[-.3em]&\label{Euler-fluid-accret-ass5}
\hspace{7.5em}\exists K<+\infty\,{:}\ \ \ 
\Big|\frac{\theta\,\psi_\theta'(J,\theta)}J-\,\psi_J'(J,\theta)\Big|
\le  K\Big(1{+}\frac{\phi(J)}J+\W(J,\theta)\Big)\,,
\\[-.4em]\nonumber
&v_{\rm ext}^{}\in C(\R^+{\times}\varOmega{\times}[0,+\infty))\,,\ \
\vv_{\rm ext}^{}\in C(\R^+{\times}\varOmega;\R^3)\,:\ \ 
\\[-.4em]&\hspace{7.5em}
|v_{\rm ext}^{}(t,\xx,J)|+v_{\rm ext}^{}(t,\xx,J)|\vv_{\rm ext}^{}(t,\xx)|^2\le K\,,\ \
\label{Euler-fluid-accret-ass-v}
\\&\label{Euler-fluid-accret-ass2}
  \W=\W(J,\theta)\ge\epsilon(1+\theta^\beta)/J\,,\ \ \ 
\DIS(J,\theta;\ee){:}\ee\ge0\,,
\ \ \ \ h_{\rm ext}\ge0\,,
\\&\label{Euler-fluid-accret-ass3}
\DIS(J,\theta;\ee){:}\ee\ge\epsilon|\ee|^q\,,\ \ \kappa(J,\theta)\ttt\ge\epsilon\eee
%G{\cdot}G\ge\epsilon|C|^2
\,,
\\&\nonumber
0<\inf\varrho_0\le\sup\varrho_0<+\infty\,,\ \ \ \
\vv_0\in L^2(\varOmega;\R^3)\,,\ \ \
0\le\theta_0\in L^1(\varOmega)\,,\ \
\\[-.2em]&\hspace{13em}
\text{and }\ \ \phi(J_0)/J_0+\W(J_0,\theta_0)\in L^1(\varOmega)\,.
\label{Euler-fluid-accret-ass4}\end{align}\end{subequations}
In \eq{Euler-fluid-accret-ass2}, $\W$ is from \eq{Euler4-fluid-accret}.
Since $\phi'$ does not occur in the system \eq{Euler-fluid-accret},
discontinuities as in Figure~\ref{PT-accretion} can be considered. 

% \COMMENT{OUR PAPER $\mbox{\cite{MiRoSt23MGDC}}$ IS NOW MODIFIED:}\DELETE{Compared with $\mbox{\cite{MiRoSt23MGDC}}$, the estimation strategy here
% follows a different path, as  we control the mass density taking advantage
% of the singular character of $\phi$ under strong compression, namely,
% \eqref{compressibility}.}
Like  in \cite{Roub24TVSE}, the strategy of the proof is first to read some formal
apriori estimates from the total-energy balance \eq{total-engr-selfgrav-accret}.
To this aim, we need to consider (and then to estimate) the negative gravitational
energy $\varrho\GRAVPOT$ on the right-hand side. \ttt By neglecting \eee the energy of the
gravitational field \ttt which is surely non-negative at each time
instant\eee, we consider the inequality
\begin{align}\nonumber 
&\int_\varOmega\frac{\varrho(t)}{2}|\vv(t)|^2
+\frac{\!\phi(J(t))\!}{J(t)}+\W(J(t),\theta(t))\,\d\xx\le
\!\int_0^t\!\!\int_\varOmega\!\bigg(r_{\rm ext}^{\nothing}(\vv_{\rm
                   ext}^{\nothing}{-}\vv){\cdot}\vv 
+\frac{r_{\rm ext}^{\nothing}}2|\vv|^2-\varrho\GRAVPOT 
\\[-.3em]&\hspace*{4em}
\ttt+r_{\rm ext}^{}\GRAVPOT\eee
+v_{\rm ext}^{}(J)\Big(\sss \frac{\psi(J,\theta)-\theta \psi_\theta'(J,\theta)
           }{J} \eee\Big) 
+h_{\rm ext}^{\nothing}\!\bigg)\d\xx\d t-\int_\varOmega\varrho(t)\GRAVPOT(t)\,\d\xx
+C_0 \,
 \label{total-engr-selfgrav-accret-modif}\end{align}
with $C_0$  depending on the initial conditions as qualified in
\eq{Euler-fluid-accret-ass4}.%  Note that the free-energy term in the
% right-hand side above reads
% $$ 
% \frac{\psi(J,\theta)}{J}-\psi_J'(J,\theta)=\frac{\phi(J)}J+\W(J,\theta)
% +\theta\,\frac{\psi_\theta'(J,\theta)}J\,-\,\psi_J'(J,\theta)
% $$
% due to \eq{Euler-fluid-accret-ass5}.

The gravitational energy in \eq{total-engr-selfgrav-accret-modif} can be
estimated by the Young inequality as
\begin{align}
-\int_\varOmega\varrho\GRAVPOT\,\d\xx\le\frac12\int_\varOmega
\varrho^2+\GRAVPOT^2\,\d\xx\,.
\label{total-engr-selfgrav-accret-modif+}\end{align}
The last term in \eq{total-engr-selfgrav-accret-modif+} is to be
estimated by using
\begin{align}(-\GRAVPOT(\xx))^r\stackrel{\mbox{\footnotesize\eq{V-integral}}}{=}
\GRAVCONST^r\bigg(\!\int_{\varOmega}
\frac{\varrho(\widetilde\xx)}{|\xx{-}\widetilde\xx|}\,\d\widetilde\xx\bigg)^{\!r}\stackrel{\text{H\"older}}{\le}
 \GRAVCONST^r\bigg(\!\!\!\!\lineunder{\sup_{\xx\in\varOmega}
\int_{\varOmega}\frac1{|\xx{-}\widetilde\xx|^{r/(r-1)}}\,\d\widetilde\xx}
{$=: C_{r,\varOmega}$}\!\!\!\bigg)^{\!r-1}\!\!\int_{\varOmega}\!\varrho(\widetilde\xx)^r\,\d\widetilde\xx\,,
  \nonumber\\[-1.7em]\label{V-integral-est}\end{align} 
where the constant $C_{r,\varOmega}$ depends on the bounded domain $\varOmega$
and on the exponent $r$ and is indeed finite if $r>3/2$, while
$C_{r,\varOmega}=+\infty$ for $r\le3/2$. We use it in
\eq{total-engr-selfgrav-accret-modif+} with $r=2$. 
Then, using also \eq{cont-eq+}, we estimate the right-hand side of
\eq{total-engr-selfgrav-accret-modif+} as
$(1+C_{2,\varOmega})\sup_{\varOmega}^{}\varrho_0\int_\varOmega1/J^2\,\d\xx$, which
eventually allows for  absorbing the  term $\phi(J)/J$  in
the left-hand side 
if $\phi$ has a growth faster than $1/J$ for $J\to0+$, as assumed in
\eq{compressibility}; cf.\ also Fig.~\ref{PT-accretion}.

The term $r_{\rm ext}^{\nothing}(\vv_{\rm ext}^{\nothing}{-}\vv){\cdot}\vv$
in \eq{total-engr-selfgrav-accret-modif} can be estimated
by using \eq{r-vs-v} as
\begin{align}
\int_\varOmega&\!r_{\rm ext}^{\nothing}(\vv_{\rm ext}^{\nothing}{-}\vv){\cdot}\vv\,\d\xx
=\int_\varOmega\!\sqrt{r_{\rm ext}^{\nothing}}\vv_{\rm ext}^{\nothing}
{\cdot}\sqrt{v_{\rm ext}^{\nothing}\varrho}\vv
-v_{\rm ext}^{\nothing}\varrho|\vv|^2\,\d\xx
\le\int_\varOmega\!r_{\rm ext}^{\nothing}|\vv_{\rm ext}^{\nothing}|^2\!
+2v_{\rm ext}^{\nothing}\varrho|\vv|^2\,\d\xx\,,
\label{est-rhs2}\end{align}
where, using also \eq{Euler-fluid-accret-ass-v}, the last term can be treated
by the Gronwall inequality while $r_{\rm ext}^{\nothing}|\vv_{\rm ext}^{\nothing}|^2$
is bounded due to \eq{Euler-fluid-accret-ass-v}. Similarly,
$\frac12r_{\rm ext}^{\nothing}|\vv|^2=\frac12v_{\rm ext}^{\nothing}\varrho|\vv|^2$,
as well as the   other right-hand side terms in
\eq{total-engr-selfgrav-accret-modif} can be treated by the Gronwall
inequality,  also using assumption \eq{Euler-fluid-accret-ass5}.
 We thus obtain  uniform-in-time bounds on
$\varrho|\vv|^2\in L^1(\varOmega)$, $\phi(J)/J\in L^1(\varOmega)$, and
$\W\in L^1(\varOmega)$ also using the fact that $\theta \geq 0$ from the heat
equation.  Here, we also used  the first condition in \eq{Euler-fluid-accret-ass2}.

Remarkably, the above estimates rely on (\ref{Euler-fluid-accret-ass}a--d) only.
Namely, they hold also in absence of viscosity.  This allows to assume that
viscosity degenerates (similarly as the elastic response) in the low-density
dilute area.

For further estimates, we still use the energy-dissipation balance
\eq{mech-engr-selfgrav-accret}.
In comparison with \eq{total-engr-selfgrav-accret-modif}, now we have
the dissipation-rate term in the left-hand side, which gives the estimate on
$\ee(\vv)$.  Further, one estimates  $v_{\rm ext}^{}(J)(\phi(J)/J-\phi'(J))$
in terms of  $\phi(J)/J$ and then  applies  the Gronwall
inequality.  Eventually,   we estimate the adiabatic right-hand side
term $(\psi_J'(J,\theta){-}\phi'(J))
\,{\rm div}\,\vv$, for which we use \eq{Euler-fluid-accret-ass3} 
together with the already obtained estimates. From \eq{Euler-fluid-accret-ass3},
we obtain also the a-priori bound on
$\ee(\vv)\in L^q([0,T]{\times}\varOmega;\R^{3\times 3})$.

Moreover, having  already proved that  the heat sources 
are in  $L^1([0,T]{\times}\varOmega)$,
we  obtain the additional  integrability of $\nabla\theta$ in
$L^r([0,T]{\times}\varOmega;\R^{3})$ for $1\le r<(3{+}2\beta)/(3{+}\beta)$
with $\beta>1$ from \eq{Euler-fluid-accret-ass2}; cf.\ \cite{Roub24TVSE}.

\begin{remark}[{\sl Analysis towards existence of solutions}]\label{rem-hyper}\upshape
 The  above estimates  indicate a possible stability frame for a numerical
 scheme, for instance by a finite-element discretization of the momentum
and the heat equations.  To prove the convergence of approximations and the
existence of weak solution one would follow the strategy in
\cite{DucFei04DGS} by extending it to an open system. This can be
expedited by considering a variant of  the model, featuring also 
some higher-gradient  terms. In particular, one could consider
a stronger viscosity, following the theory by E.\,Fried and
M.\,Gurtin \cite{FriGur06TBBC}, as already considered in the general
nonlinear context of {\it multipolar fluids} by J.\,Ne\v cas at al.\ 
\cite{Neca94TMF,NeNoSi89GSIC,NecSil91MVF} and as originally inspired by
R.A.\,Toupin \cite{Toup62EMCS} and R.D.\,Mindlin \cite{Mind64MSLE}.
More specifically, the local dissipation potential $\DIS(J,\theta;\cdot)$
may be enhanced by a nonlocal term as
$$
\mathscr{D}(J,\theta;\ee(\vv),\nabla^2\vv)
=\int_\varOmega\DIS(J,\theta;\ee(\vv))+\frac1q\nu(J,\theta)|\nabla^2\vv|^q\,\d\xx
$$
with $q>3$ and $\nu(\cdot,\cdot)>0$ continuous with positive infimum
on $\R{\times}\R$. This brings an additional symmetric hyper-stress contribution
 to the Cauchy stress, namely ${\rm div}(\nu|\nabla^2\vv|^{q-2}\nabla^2\vv)$.
 The dissipation rate $\xi$ in \eq{entropy-eq} and also the dissipation rate in
 the energy-dissipation balance \eq{mech-engr-selfgrav-accret} expands by
 $\nu|\nabla^2\vv|^p$,  and the  boundary conditions should be
 extended appropriately. Note that such  hyper-viscosity may be seen as not
 entirely physical in  the case of a very diluted system, regardless of the
 size of $\nu$. On the other hand, the above derived a-priori
 estimates are uniform with respect to $\nu\to0$,  so that the
 addition of the second-order contribution is compatible with 
 the global energetic of the system.
 Refer to \cite{Roub24TVSE} for  the technically nontrivial analysis 
 for a simpler,  albeit  similar thermomechanical system.
\end{remark}

\begin{example}\upshape
 An  example of $\psi$  fulfilling  the assumptions
(\ref{Euler-fluid-accret-ass}a,b,d)   with  $\phi=\psi(\cdot,0)$
 as in  Figure~\ref{PT-accretion}  is 
\begin{align}\label{exa}
\psi(J,\theta)=\max\Big(0,\frac1{J^\alpha}-1\Big)+b\frac\theta{J^z}
+c_0\theta\big(1{-}{\rm ln}\,\theta\big)-c_1\theta^{\beta}
\end{align}
with some $\alpha>1$, $\beta>1$,  $z >0$, $c_0>0$, and $c_1>0$. 
Note that $\W(J,\theta)=c_0\theta/J+c_1(\beta{-}1)\theta^{\beta}/J
$ so that also \eq{Euler-fluid-accret-ass2} is satisfied.

The heat capacity reads 
$c(J,\theta)=c_0/J+c_1\beta(\beta{-}1)\theta^{\beta-1}/J$.  In
particular, the term  $-b\theta/J^z$  does not contribute to
the heat capacity, although it provides adiabatic heating during compression.

 In order to check \eq{Euler-fluid-accret-ass5},  we compute
\begin{align*}
  &\frac{\theta\psi_\theta'(J,\theta)}{J}-\psi_J'(J,\theta) =-\frac{(z{+}1)b\theta}{J^{z+1}}
-\frac{c_0\theta{\rm ln}\theta}{J}-\frac{\beta c_1\theta^{\beta}}{J} +
\left\{
  \begin{array}{ll}
   \displaystyle \frac{\alpha}{ J^{\alpha+1}}\quad&\text{for}
                                                           \ \ J\leq 1,\\
    0&\text{for} \ \ J>1.
  \end{array}
  \right. 
\end{align*}
 This can indeed be
dominated by $\phi(J)/J+\W(J,\theta)$ as required in \
\eq{Euler-fluid-accret-ass5}. 

 In fact, controlling the pressure $-\psi_J'(J,\theta)$ by the
actual stored energy $\phi(J)/J$ (or, equivalently,  controlling  the
{\it Kirchhoff} pressure $-J\psi_J'(J,\theta)$ by the referential stored
energy $\phi(J)$) corresponds to the control of the Kirchhoff stress
$\varphi'(\FF)\FF^\top=-J\psi_J'(J,\theta)\bbI$ by the referential stored
energy $\varphi(\cdot,0)$. This may be achieved via an
application of the Young inequality   
$$ \frac{\theta}{J^{z}}\le\beta \theta^\beta+\frac{\beta-1}{\beta}\frac{1}{J^{z\beta/(\beta-1)}},$$  and calls for asking  $\alpha\ge z\beta/(\beta-1)$. The
relevance of the possibility of controlling the Kirchhoff stress via the
stored energy was pointed out by J.M.\,Ball \cite{Ball84MELE,Ball02SOPE}.
Some analogous argument applies to  the control  of the term
$\theta\psi_\theta'(J,\theta)=-b\theta/J^{z}
-c_0\theta{\rm ln}\theta-\beta c_1\theta^{\beta}$.
\end{example}

\begin{remark}[{\sl Existence of the deformation $\yy$}]\label{rem-invert}\upshape
   Let us close this discussion by commenting on the possibility of recovering a
  deformation $\yy(t,\cdot)$ from our modelization in actual coordinates, i.e.,
  from $\bm\xi(t,\cdot)$. Note at first that   the blow-up of
$\phi(\cdot,\theta)$ at $0$ ensures that $J=\det\FF>0$, so that the
local invertibility of $\bm\xi(t,\cdot)$ ensues. In order
$\bm\xi(t,\cdot)$ to be globally invertible, some stronger assumption would be
required. In particular, by assuming the stronger boundary condition $\vv ={\bm
  0}$ on $\varGamma$, one would have that $\bm\xi(t,\xx)=\xx$ at
$\varGamma$. Some stronger coercivity of the elastic energy would then
allow to follow the approach in \cite{Ball81GISF} and prove the global
invertibility of $\bm\xi(t,\cdot)$. The same would hold as long as
$\bm\xi(t,\cdot)$ on $\varGamma$ coincides with the trace of a
homeomorphism, possibly being different from the identity. 
 \end{remark}

\bigskip\bigskip

{\small
\baselineskip=12pt

\noindent{\it Acknowledgments.}
Fruitful discussions with \v{S}\'arka Ne\v{c}asov\'a \ttt and a
very careful reading with useful suggestions and correction of some
formulae by an anonymous referee \eee are greatfully
acknowledged. A support from the CSF project GA23-04676J and the institutional
support RVO: 61388998 (\v CR) are gratefully acknowledged. US is supported
by the Austrian Science Fund (FWF) through projects 10.55776/F65,
10.55776/I4354, 10.55776/I5149, and 10.55776/P3278.

\bibliographystyle{plain}

} 
\end{document}